\documentclass[journal,10pt,onecolumn,draftclsnofoot,]{IEEEtran}
\ifCLASSINFOpdf
\else
\fi
\usepackage{graphicx}          
\usepackage{psfrag}
\usepackage{enumerate}
\usepackage{overpic}
\usepackage{euscript}
\usepackage{graphicx}
\usepackage{amssymb}
\usepackage{amsmath}
\usepackage{algorithm}
\usepackage{psfrag}
\usepackage[noend]{algpseudocode}
\usepackage{color,comment}
 \newtheorem{defn}{Definition}
 \newtheorem{rem}{Remark}
\newtheorem{thm}{Theorem}
\newtheorem{lem}{Lemma}
\newtheorem{assumption}{Assumption}
\def\R{\mathbb{R}}

\def\Wb{\mathbf{W}}
\def\Z{\mathbb{Z}}
\def\bmu{\boldsymbol{\mu}}

\def\bx{\boldsymbol{x}}
\def\bbx{\bar{\boldsymbol x}}
\def\by{\boldsymbol{y}}
\def\bba{\mathbf{a}}

\def\ba{\boldsymbol{a}}

\def\Qb{\boldsymbol{Q}}
\def\Rb{\boldsymbol{R}}
\def\Fb{\boldsymbol{F}}
\def\fb{\boldsymbol{f}}
\def\zq{\boldsymbol{q}}
\def\Mb{\boldsymbol{M}}
\def\xib{\boldsymbol{\xi}}
\def\ra{\textrm{a}}
\def\s{~}
\newcommand{\E}{\mathbb{E}}
\newcommand{\Tr}{\mathrm{Tr}}

\def\zbx{\mathbf{X}}
\def\ze{\boldsymbol e}
\definecolor{amber}{rgb}{1.0, 0.49, 0.0}
\def\revi#1{{\color{black}#1}}

\def\re1#1{{\color{black}#1}}
\def\r2#1{{\color{black}#1}}
\def\rev3#1{{\color{black}#1}}

\begin{document}

\title{Stochastic Learning in Potential Games: Communication and Payoff-Based Approaches}

\author{Tatiana~Tatarenko
\thanks{T. Tatarenko is with the Control Methods and Robotics Lab 
Technical University Darmstadt, Darmstadt, Germany}}

\maketitle
\begin{abstract}                          
Game theory serves as a powerful tool for distributed optimization in
multi-agent systems in different applications. In this paper we consider multi-agent systems that can be modeled by means of potential games whose
potential function coincides with a global objective function to be maximized.
In this approach, the agents correspond to the strategic decision makers and the optimization problem is equivalent to the problem
of learning a \rev3{potential function maximizer} in the designed game. The paper deals with two different information settings in the system. Firstly, we consider systems, where agents have the access to the gradient of their utility functions. However, they do not possess the full information about the joint actions. Thus, to be able to move along the gradient toward a local optimum, they need to exchange the information with their neighbors by means of \emph{communication}. The second setting refers to a \emph{payoff-based} approach to learning potential function maximizers. Here, we assume that at each iteration agents can only observe their own played actions and experienced payoffs. \revi{In both cases, the paper studies unconstrained non-convex optimization with a differentiable objective function.}
To develop the corresponding algorithms guaranteeing convergence to a \rev3{local maximum of the potential function} in the game, we utilize the idea of the well-known Robbins-Monro procedure based on the theory of stochastic approximation.
\end{abstract}

\section{Introduction}
The goal of distributed optimization in a multi-agent system consists in establishment of the agents' local dynamics leading to a common goal.
\revi{This paper deals with unconstrained optimization problems of the following type:
 \begin{align}\label{eq:opt0}
\phi(\ba) = \phi(a_1,\ldots,a_N)\to\max,\quad \ba\in\R^N,
 \end{align}
 where $N$ is the number of agents in a system and $a_i \in \R$\footnote{One-dimensional case is chosen for the sake of notational simplicity. All results in this paper can be extended to any finite dimensional case.} is the variable in the objective function $\phi$ controlled by the agent $i$.}
 A relation between the optimization problem above and \emph{potential games} has been recently established \cite{MardenMist}, \cite{MardDes}, \cite{MardOver}, \cite{Olfati}.
This allows for using potential games to find a solution to the problem \eqref{eq:opt0}. Much work on modeling a multi-agent system as a potential game has been carried out \cite{MardOver}, \cite{StrLearn}.
The proposed approach is shown to be applicable to such optimization problems as consensus finding \cite{MardenMist}, sensor coverage \cite{Goto_PIPIP}, \cite{COVEr},
vehicle-target assignment \cite{Vehicle}, broadcast tree formation \cite{Tardos}, routing over networks \cite{rosenthal}, \cite{Roughgarden}, just to name a few.

The idea of using potential games in distributed multi-agent optimization is the following. Firstly, the local utility functions of agents need to be modeled in such a way that the maximizers of the global objective function coincide with potential function maximizers in the resulting game. After that one needs to develop a learning procedure, which leads agents to one of these states.
Among learning algorithms that are applicable to potential games and have been presented in the literature so far the following ones demonstrate an efficient performance by approaching the system to its optimum as time runs:
the \emph{log-linear learning} \cite{MardRev}, \cite{ACC}, its payoff-based and synchronous versions \cite{MardRev}, \emph{payoff-based inhomogeneous partially irrational play} \cite{Goto_PIPIP}, \cite{COVEr}, and \emph{adaptive $Q$-learning and $\varepsilon$-greedy decision rule} \cite{Chapman}. All these algorithms can be applied only to discrete optimization problems and, thus, assume players to have discrete actions. This fact restricts the applications of the algorithms above in control, engineering, and economics, \r2{since in these areas there are a lot of optimization problems formulated in continuous settings due to the nature of control, strategy, or price variables (for numerous examples, see \cite{A1}, \cite{Perkins} and references
therein).}
On the other hand, some learning procedures have been proposed in \cite{Benaim2015}, \cite{Krichene2015}, \cite{MardDes}, \cite{Perkins}, \cite{Ratliff2013} that allow us to deal with systems with continuous states. However, the efficiency of the \emph{gradient play} \cite{MardDes} is strictly connected to the communication
abilities of agents, their access to the closed form of the utilities' gradients, and \emph{convexity} of the objective function. The \emph{mixed-strategy learning} procedures presented in \cite{Benaim2015}, \cite{Krichene2015}, \cite{Perkins}, in its turn, require each agent to have access to a so called \emph{oracle information}, namely to be aware about her payoff that she would have got for each action against the selected joint action of others. The access to this information might be provided in routing games, but it can be still questionable in many other applications. Some methods for tackling continuous games have been presented in \cite{Ratliff2013}. However, all of them use the information on second order derivatives.

In this paper we propose two learning algorithms, namely \emph{communication-based} and \emph{payoff-based} algorithms, to learn local maxima of potential functions in potential games.
{In contrast to the works mentioned above, we deal with \emph{continuous action potential games}, consider the case of \emph{non-concave potential functions} \revi{with bounded Lipschitz continuous gradients}, and assume that \emph{neither oracle information nor second order derivatives} are available for agents.}
In the approach based on communication we assume agents to be connected with each other by means of time-varying communication graph. By utilizing this graph they exchange the information with their neighbors and, thus, update their estimations of the current joint action. It allows them to use these estimations to update their actions and move toward the gradient ascent. As we prefer not to introduce the assumption on double stochasticity of the communication graph and, thus, not to complicate the algorithm's implementation \cite{NonDS}, the \emph{push-sum protocol} \cite{16A2}, \cite{Tsianos2012}, which is based only on the nodes' out-degrees, is adopted to learn local optima in potential game.
Optimization problems that require a communication-based approach often arise in signal processing applications in sensor networks \cite{LiHan2010}, \cite{A1}, \cite{Neglia2009}.
However, if communication is not set up in the system and the agents have no access to the gradient information, a payoff-based algorithm needs to be applied to the formulated optimization problem. In the payoff-based approach agents can only observe their currently played actions and the corresponding utility values. Based on this information agents iteratively update their mixed strategies. To model the mixed strategies of the agents, we use the Gaussian distributions over the action sets, following the idea presented in the literature on the learning automata \cite{Beigy}, \cite{Thatha}. Payoff-based learning is of interest in various applications encountering uncertainties of the physical world \cite{Hatanaka2016} or possessing a complex problem formulation without any closed-form expression \cite{Windfarm}.

Since we \emph{do not assume concavity of the potential function}, the main goal we pursue in both algorithms is to escape from those critical points of the global objective function which are not local maxima.
\revi{Many papers on distributed non-convex optimization prove only convergence to some critical point \cite{BianchiStochApp}, \cite{NonDS}, \cite{NCScutari}.
Some works in centralized non-convex optimization study escaping from saddle points \cite{nonconvexGE}, \cite{Anandkumar2016EfficientAF}, \cite{pemantle1990}.
However, as it is shown in \cite{Anandkumar2016EfficientAF}, the problem of escaping saddle points of the order four or higher is NP-hard, while for escaping from saddle points of a lower order existence of the second or even the third order derivatives of the objective function is required.
This work aims to propose algorithms that enable escaping from local minima without the assumption on existence of the second order derivatives.}
To reach this goal we add a \emph{stochastic} term to the iterative learning steps in the communication-based algorithm (the payoff-based algorithm is stochastic itself) and refer to the results on the Robbins-Monro procedure, which is extensively studied in the literature on stochastic approximation, particularly in \cite{NH}. Two main results are formulated in this paper: the first one claims the almost sure convergence of the communication-based algorithm to a critical point \revi{different from local minima}, whereas the second one claims the convergence in probability of the payoff-based algorithm to such a state.

The paper is organized as follows. Section\s\ref{sec:prel} provides some preliminaries that are required further for the theoretic analysis. Section\s\ref{sec:gta} formulates a game-theoretic approach in multi-agent optimization and addresses information restrictions.
Section\s\ref{subsec:comba} introduces and analyzes the communication-based algorithm for \revi{unconstrained} optimization problems modeled by means of potential games with continuous actions. Section\s\ref{sec:pba} deals with the payoff-based approach to such problems. Section\s\ref{sec:impl} contains an illustrative implementation of the proposed algorithms. Finally, Section\s\ref{sec:concl}  draws the conclusion.

\textbf{Notation:} We will use the following notations throughout this manuscript: We denote the set of integers by $\Z$ and the set of non-negative integers by $\Z^+$. For the metric $\rho$ of a metric space $(X,\rho(\cdot))$ and a subset $B\subset X$, we let $\rho(x,B)=\inf_{y\in B}\rho(x,y)$. We denote the set $\{1,\ldots,N\}$ by $[N]$. We use boldface to distinguish between the vectors in a multi-dimensional space and scalars. Unless otherwise specified, $\|\cdot\|$ denotes the standard Euclidean norm in some space $\R^d$. The dot product of two vectors $\bx, \by\in\R^d$ is denoted by $(\bx,\by)$.  Throughout this work, all the time indices such as $t$ belong to $\Z^+$.

\section{Preliminaries}\label{sec:prel}
\subsection{Stochastic Approximation}
In the following we will use the results on the convergence of stochastic approximation Robbins-Monro procedure presented in \cite{NH}.
We start by introducing some important notations.
Let $\{\zbx(t)\}_t$, $t\in \Z^+$, be a Markov process on some state space $E\subseteq \R^d$. The transition function of this process, namely $\Pr\{\zbx(t+1)\in G| \zbx(t)=\zbx\}$, is denoted by $P(t,\zbx,t+1,G)$, $G\subseteq E$.

\begin{defn}\label{def:def1}
The operator $L$ defined on the set of measurable functions $V:\Z^+\times E\to \R$, $\zbx\in E$, by
\begin{align*}
LV(t,&\zbx)=\int{P(t,\zbx,t+1,dy)[V(t+1,y)-V(t,\zbx)]}\cr
&=E[V(t+1,\zbx(t+1))\mid \zbx(t)=\zbx]-V(t,\zbx).
\end{align*}
is called a \emph{generating operator} of a Markov process $\{\zbx(t)\}_t$.
\end{defn}

Let $B$ be a subset of $E$, $\EuScript U_{\epsilon}(B)$ be its $\epsilon$-neighborhood, i.e. $\EuScript U_{\epsilon}(B)=\{\zbx: \rho(\zbx,B)<\epsilon\}$.
Let $\EuScript W_{\epsilon}(B)=E\setminus \EuScript U_{\epsilon}(B)$ and $\EuScript U_{\epsilon,R}(B)=\EuScript U_{\epsilon}(B)\cap \{\zbx:\|\zbx\|<R\}$.

\begin{defn}\label{def:def2}
The function $\psi(t,\zbx)$ is said to belong to class $\Psi(B)$, $\psi(t,\zbx)\in\Psi(B)$, if

1) $\psi(t,\zbx)\ge 0$ for all $t\in \Z^+$ and $\zbx\in E$,
\smallskip

2) for all $R>\epsilon >0$ there exists some $Q=Q(\epsilon, R)$ such that
$\inf_{t\ge Q, \zbx\in \EuScript U_{\epsilon,R}(B)}\psi(t,\zbx)>0$.
\end{defn}

Further we consider the following general recursive process $\{\zbx(t)\}_t$ taking values in $\R^d$:
\begin{align}\label{eq:basic}
 \zbx(t+1) = \zbx(t) & + \alpha(t+1){\Fb}(t,\zbx(t)) \cr
 &  + \beta(t+1)\Wb(t+1,\zbx(t),\omega),
\end{align}
where $\zbx(0)=\bx_0\in \R^d$, ${\Fb}(t,\zbx(t))=\fb(\zbx(t))+\zq(t,\zbx(t))$ such that
$\fb: \mathbb{R}^d\to\mathbb{R}^d$, $\zq(t,\zbx(t)): \mathbb{Z}^{+}\times\mathbb{R}^d\to\mathbb{R}^d$,
$\{\Wb(t,\bx,\omega)\}_t$ is a sequence of random vectors such that
\begin{align*}
\Pr \{\Wb(t+1,\zbx(t),\omega)&|\zbx(t),\ldots,\zbx(0)\}\cr
&=\Pr\{\Wb(t+1,\zbx(t),\omega)|\zbx(t)\},
\end{align*}
 and $\alpha(t)$, $\beta(t)$ are some positive step-size parameters. Under the assumption on the random vectors $\Wb(t)$ above, the process \eqref{eq:basic} is a Markov process.
We use the notation $A(t,\bx)$ to denote the matrix with elements $\E W_i(t+1,\bx,\omega) W_j(t+1,\bx,\omega)$.

Now we quote the following theorems for the process \eqref{eq:basic}, which are proven in \cite{NH} (Theorem 2.5.2 and Theorem 2.7.3 respectively).

\begin{thm}\label{th:finiteness}
  Consider the Markov process defined by \eqref{eq:basic} and suppose that there exists a function $V(t,\zbx)\ge 0$ such that $\inf_{t\ge0}V(t,\zbx)\to\infty$ as $\|\zbx\|\to\infty$ and
  \[LV(t,\zbx)\le -\alpha(t+1)\psi(t,\zbx) + g(t)(1+V(t,\zbx)),\]
   where $\psi\in\Psi(B)$ for some set $B\subset\R^d$, $g(t)>0$, $\sum_{t=0}^{\infty}g(t)<\infty$. Let $\alpha(t)$ be such that $\alpha(t)>0$, $\sum_{t=0}^{\infty} \alpha(t)= \infty$.
   Then almost surely $\sup_{t\ge 0}\|\zbx(t)\| = R< \infty$.
\end{thm}

\begin{thm}\label{th:th1}
  Consider the Markov process defined by\s\eqref{eq:basic} and suppose that $B=\{\zbx: \fb(\zbx)=\boldsymbol 0\}$ be a union of finitely many connected components. Let the set $B$, the sequence $\gamma(t)$, and some function $V(t,\zbx): \mathbb{Z}^+\times\mathbb{R}^d\to\mathbb R$ satisfy the following assumptions:

 \begin{enumerate}[(1)]
 	\item\label{th1:cond1}  For all $t\in \Z^+$, $\zbx\in \mathbb{R}^d$, we have $V(t,\zbx)\ge 0 $ and $\inf_{t\ge0}V(t,\zbx)\to\infty$ as $\|\zbx\|\to\infty$,
 	\medskip
 	
	\item\label{th1:cond2}  $LV(t,\zbx)\le -\alpha(t+1)\psi(t,\zbx) + g(t)(1+V(t,\zbx))$, where  $\psi\in\Psi(B)$, $g(t)>0$, $\sum_{t=0}^{\infty}g(t)<\infty$,
	\medskip
	
	\item\label{th1:cond3}  $\alpha(t)>0$, $\sum_{t=0}^{\infty} \alpha(t)= \infty$, and $\sum_{t=0}^{\infty} \alpha^2(t)<\infty$, and there exists $c(t)$ such that $\|\zq(t,\zbx(t))\|\le c(t)$ almost surely (a.s.) for all $t\ge0$ and $\sum_{t=0}^{\infty}\alpha(t+1)c(t)<\infty$,
	\medskip
	
	\item\label{th1:cond4} $\sup_{t\ge0,\|\zbx\|\le r}\|{\Fb}(t,\zbx)\|<\infty$ for any $r$, and
	\medskip
	
	\item\label{th1:cond5} $\E \Wb(t) = \boldsymbol 0$, $\sum_{t=0}^{\infty}\E\beta^2 (t)\|\Wb(t)\|^2<\infty$.
\end{enumerate}
{Then the process~\eqref{eq:basic} converges almost surely either to a point from $B$ or to the boundary of one of its connected components, given any initial state $\zbx(0)$}.
\end{thm}

Thus, if $\fb(\bx)$ is the gradient of some objective function $\Phi(\bx)$: $\R^d\to\R$, then the procedure \eqref{eq:basic} converges almost surely to a critical point of $\Phi$ as $t\to\infty$.
However, the theorem above cannot guarantee the almost sure convergence to a local maximum of the objective function $\Phi$. To exclude the convergence to critical points that are not local maxima, one can use the following theorem proven in \cite{NH} (Theorem 5.4.1).

\begin{thm}\label{th:th5}
Consider the Markov process $\{\zbx(t)\}_t$ on $\mathbb{R}^d$ defined by \eqref{eq:basic}.
Let $H'$ be the set of the points $\bx'\in\mathbb{R}^d$ for which there exists a symmetric positive definite matrix $C=C(\bx')$  and a positive number  $\epsilon=\epsilon(\bx')$ such that $(\fb(\bx),C(\bx-\bx'))\ge 0$  for $\bx\in U_{\epsilon}(\bx')$.
Assume that for any $\bx{'}\in H{'}$ there exists positive constants $\delta=\delta(\bx')$ and $K=K(\bx')$ such that
\begin{enumerate}[(1)]
 \item\label{th5:cond1} $\|\fb(\bx)\|^2 + |\Tr[A(t,\bx)-A(t,\bx')]|\le K\|\bx - \bx'\|$
 for any $\bx:$ $\|\bx - \bx'\|<\delta$,
 \medskip

 \item\label{th5:cond2} $\sup_{t\ge 0, \|\bx - \bx'\|<\delta}\E\|\Wb(t,\bx,\omega)\|^4<\infty$.
  \medskip
\end{enumerate}
  Moreover,

\label{th5:cond3}{(3)} $\sum_{t=0}^{\infty}\beta^2(t)<\infty$, $\sum_{t=0}^{\infty}\left(\frac{\beta(t)}{\sqrt{\sum_{k=t+1}^{\infty}\beta^2(k)}}\right)^3<\infty$,
  \medskip

\label{th5:cond4}{(4)} $\sum_{t=0}^{\infty}\frac{\alpha(t)q(t)}{\sqrt{\sum_{k=t+1}^{\infty}\beta^2(k)}}<\infty$, $q(t)=\sup_{\bx\in\mathbb{R}^d}\|\zq(t,\bx)\|$.

Then $\Pr\{\lim_{t\to\infty}\zbx(t)\in H'\}=0$ irrespective of the initial state $\zbx(0)$.
\end{thm}

\subsection{Consensus under Push-Sum Protocol}\label{subsec:pertps}

Here we discuss the general push-sum algorithm initially proposed in \cite{16A2} to solve consensus problems, applied in \cite{A1} to distributed optimization of convex functions, and analyzed in \cite{TAC} in context of non-convex distributed optimization. The communication-based algorithm that will be introduced in Section\s\ref{subsec:comba} is its special case.
\revi{However, to be able to apply the push-sum algorithm to the problem \eqref{eq:opt} we will need to adjust the iteration settings as well as the analysis of the procedure's behavior.}
Consider a network of $N$ agents. At each time $t$, node $i$ can only communicate to  its  out-neighbors  in  some  directed  graph $G(t)$, where the graph $G(t)$ has the vertex  set $[N]$ and the edge set $E(t)$. We introduce the following standard definition for the sequence $G(t)$.
\begin{defn}\label{def:CG}
	The sequence of communication graphs $\{G(t)\}$ is $S$-strongly connected, i.e. for any time $t\geq 0$, the graph
	\[G(t:t+S)=([N],E(t)\cup E(t+1)\cup \cdots \cup E(t+S-1)),\]
	is strongly connected. In other words, the union of the graphs over every $S$ time intervals is strongly connected.
\end{defn}
The $S$-strongly connected sequence of communication graphs ensures enough mixing of information among the agents \cite{A1}.
\revi{We use $N^{in}_i(t) = \{j\in[N]:\mbox{ }(j,i)\in E(t)\}$ and $N^{out}_i(t)= \{j\in[N]:\mbox{ }(i,j)\in E(t)\}$ to denote the sets of in- and out-neighbors of the agent $i$ at time $t$.} Each agent $i$ is represented by the corresponding node in $G(t)$ and is always considered to be an in- and out-neighbor of itself. We
use $d_i(t)$ to denote the out-degree of the agent
$i$, \revi{namely $d_i(t) =|N^{out}_i(t)|$}, and we assume that every agent knows its out-degree at every time $t$.
According to the general push-sum protocol, at every moment of time $t\in \Z^+$ each node $i$ maintains  vector   variables $\hat{\bba}_i(t)$, $\bx_i(t)$, $\boldsymbol w_i(t)\in \mathbb{R}^d$, as  well  as  a  scalar  variable $y_i(t):$ $y_i(0)=1$. These quantities are updated according to the following rules:
\begin{align}\label{eq:gps}
\boldsymbol w_i(t+1)&=\sum_{j\in N^{in}_i(t)}\frac{\bx_j(t)}{d_j(t)},\cr
y_i(t+1)&=\sum_{j\in N^{in}_i(t)}\frac{y_j(t)}{d_j(t)},\cr
\hat{\bba}_i(t+1)&=\frac{\boldsymbol w_i(t+1)}{y_i(t+1)},\cr
\bx_i(t+1)&=\boldsymbol w_i(t+1)+\ze_i(t+1),
\end{align}
where $\ze_i(t)$ is some $d$-dimensional, possibly random, perturbation at time $t$. As we can see, the auxiliary variables $\boldsymbol w_i(t+1)$ and $y_i(t+1)$ accumulate for each agent $i$ the state information broadcast by her in-neighbors at the time step $t$. The corresponding adjustment $\frac{\boldsymbol w_i(t+1)}{y_i(t+1)}$ is done to update $\hat{\bba}_i(t+1)$. This variable $\hat{\bba}_i(t+1)$ is introduced to cancel out any influence imbalances in the communication graph. The auxiliary state variable $\bx_i(t+1)$ in its turn can lead the estimation $\boldsymbol w_i(t+1)$ along some direction, given an adjusted perturbation $\ze_i(t+1)$. As the result of these updates, under some assumption on the perturbation terms, all individual estimations $\hat{\bba}_i(t+1)$, $i\in[N]$, reach a consensus that is equal to the running state average $\bar{\bx}(t) = \frac1N\sum_{j=1}^N \bx_j(t)$ and, thus, depends on the directions defined by these perturbation terms. Formally, the following important result is
proved in \cite{A1}.

%


\begin{thm}[\cite{A1}]\label{th:th2}
Consider  the  sequences $\{\hat{\bba}_i(t)\}_t$, $i\in [N]$, generated  by the push-sum protocol \eqref{eq:gps}. Assume that the graph sequence
$G(t)$ is $S$-strongly-connected.

(a) Let $\delta$ and $\lambda$ be the constants such that $\delta\ge\frac{1}{N^{NS}}$, $\lambda\le \left(1-\frac{1}{N^{NS}}\right)^{\frac1S}$. Then we have that a.s. for all $i$ and all $t$
\begin{align*}
 &\left\|\hat{\bba}_i(t+1)-\bar{\bx}(t) \right\| \cr
 &\le \frac{8}{\delta}\left(\lambda^t\sum_{j=1}^{N}\|\bx_j(0)\|_1+\sum_{s=1}^{t}\lambda^{t-s}\sum_{j=1}^{N}\|\ze_j(s)\|_1\right).
\end{align*}

(b) If $\{\gamma(t)\}$ is a non-increasing positive scalar sequence with $\sum_{t=1}^{\infty}\gamma(t)\|\ze_i(t)\|_1<\infty$ a.s. for all $i\in [n],$ then
$\sum_{t=0}^{\infty}\gamma(t+1)\left\|\hat{\bba}_i(t+1)-\bar{\bx}(t) \right\|<\infty$ a.s. for all $i,$
where $\|\cdot\|_1$ is the $l^1$-norm.
\end{thm}

\begin{rem}\label{rem:psprot}
 Note that the theorem above implies that if
$\lim_{t\to\infty}\|\ze(t)\|_1=0 \mbox{ a.s.},$
then $$\lim_{t\to\infty}\left\|\hat{\bba}_i(t+1)-\bar{\bx}(t) \right\|=0 \mbox{ a.s. for all } i.$$
In other words, given vanishing perturbations $\ze_i(t)$, $i\in[N]$, one can guarantee that in the push-sum algorithm all $\hat{\bba}_i(t)$ track the average state $\bar{\bx}(t)$.
\end{rem}

\section{Game-Theoretic Approach to Optimization}\label{sec:gta}
\subsection{Potential Game Design}\label{subsec:pgd}
We consider a multi-agent system consisting of $N$ agents. Each agent $i$ has her action set $A_i$ and the set of joint actions is denoted by $\boldsymbol A$. Let us suppose that the objective in the multi-agent system under consideration is captured by the function $\phi: \boldsymbol A\to\mathbb R$
such that the maxima of $\phi$ are the system optimal states. Formally, the global objective is to solve the following optimization problem (see also \eqref{eq:opt0}):
\begin{align}\label{eq:opt}
 \max \phi(\ba), \cr
 \mbox{s.t. }\ba\in \boldsymbol A.
\end{align}
Throughout the paper we assume existence of agents' local utility functions $U_i$, $i\in[N]$, ``aligned" with the global objective function $\phi$ such that the
resulting game $\Gamma(N,\boldsymbol A, \{U_i\}, \phi)$ is potential according to the following definition:
\begin{defn}\label{def:potgame}
A game $\Gamma(N,\boldsymbol A, \{U_i\}, \phi)$ is said to be \emph{potential game},
with a potential function $\phi: \boldsymbol A\to\mathbb R$, if the following condition is fulfilled:
\begin{align}\label{eq:pg}
 U_i(a_i,a_{-i})-U_i(a'_i,a_{-i})=\phi(a_i,a_{-i})-\phi(a'_i,a_{-i}),
\end{align}
for any $i\in[N]$, $a_i,a'_i\in A_i$, $a_{-i}\in A_1\times\ldots\times A_{i-1}\times A_{i+1}\times\ldots\times A_N$.
\end{defn}
Thus, the properties of system behavior can be analyzed by means of potential games \cite{Monderer}.
Note that in the case of the differentiable utility functions $U_i$, $i\in[N]$, the condition \eqref{eq:pg} is equivalent to the following one:
\begin{align}\label{eq:pg1}
 \frac{\partial \phi(\ba)}{\partial a_i}=\frac{\partial U_i(\ba)}{\partial a_i},
\end{align}
for any $i\in[N]$ and $\ba\in\boldsymbol A$.

Furthermore, we make the following assumptions regarding the problem \eqref{eq:opt}:

\begin{assumption}\label{assum:Ufgradient}
The gradient $\nabla \phi$ of the function $\phi$ exists and is bounded on $\R^N$.
\end{assumption}
Note that existence of bounded (sub)gradients is a standard assumption in optimization literature\footnote{In the case of non-differentiable functions, existence of a bounded subgradient is assumed \cite{Benaim2006}  \cite{A1},  \cite{Nesterov}.}.
\begin{assumption}\label{assum:A0}
The problem\s\eqref{eq:opt} has a solution with a finite Euclidean norm. Moreover,
the set of critical points of the objective function $\phi(\ba)$, i.e. the set of points $\ba$ such that $\nabla \phi(\ba)=\boldsymbol 0$, is a union of finitely many connected components.
\end{assumption}

Let $\boldsymbol A_0$ and $\boldsymbol A^*$ denote the set of critical points and local maxima of the function $\phi$ on $\boldsymbol A$ respectively.
Some approaches such as wonderful life utility \cite{Wolpert_WLU} and Shapley value utility \cite{SVU} have been proposed to encode the system's objective $\phi$ into agent-specific utility functions $\{U_i\}$ satisfying Definition\s\ref{def:potgame} and as result to deal with a potential game.
As it has been mentioned above, in this paper we assume existence of such utility functions and focus on appropriate algorithms that learn the set $\boldsymbol A^*$ in the corresponding potential game by utilizing only some restricted information which is described in the following subsection.

\subsection{Information in System and General Assumptions}
In the following we are focused on \emph{unconstrained optimization problems} in multi-agent systems. Thus, we deal with continuous action potential games, where $A_i=\mathbb R$, modeled for such problems according to Subsection\s\ref{subsec:pgd}. Our goal is to develop learning algorithms for two different types of information, namely communication-based and payoff-based information, available in the multi-agent system under consideration.
Before describing the information available in the system and formulating the learning algorithms, we introduce the following assumptions on the potential function in the modeled potential game under consideration.

\begin{assumption}\label{assum:Ufinfin}
 The function $-\phi$ is coercive, i.e. $\phi(\ba)\to-\infty$ as $\|\ba\|\to\infty$.
\end{assumption}

\begin{rem}\label{rem:linbeh}
  Since Assumption~\ref{assum:Ufgradient} requires the boundedness of $\nabla\phi$, the function $-\phi(\ba)$ grows not faster than a linear function as $\|\ba\|\to\infty$, given fulfillment of Assumptions~\ref{assum:Ufgradient} and\s\ref{assum:Ufinfin}. Moreover, in this case, one can assume existence of some constant $K>0$ such that $a^i\frac{\partial \phi(\ba)}{\partial a^i}\le 0$ for any $\ba:$ $|a^i|>K$.
 \end{rem}


\begin{assumption}\label{assum:A5A}
 For any point $\ba'\in \boldsymbol A_0$ that is not a local maximum of $\phi$ there exists a symmetric positive definite matrix $C(\ba')$ such that $(\nabla\phi(\ba),C(\ba')(\ba-\ba'))\le 0$ for any $\ba\in U(\ba')$, where $U(\ba')$ is some open neighborhood of $\ba'$.
\end{assumption}

\begin{rem}\label{rem:LMassum}
Note that if the second derivatives of the function $\phi$ exist, the assumption above holds for any critical point of $\phi$ that is a local minimum. Indeed, in this case
\begin{align*}
 \nabla \phi(\ba) = H\{\phi(\ba')\}(\ba - \ba')+\delta(\|\ba - \ba'\|),
\end{align*}
where $\delta(\|\ba - \ba'\|) = o(1)$ as $\ba \to \ba'$ and $H\{\phi(\cdot)\}$ is the Hessian matrix of $\phi$ at the corresponding point.
\end{rem}
\revi{However, Assumption~\ref{assum:A5A} does not require existence of second derivatives. For example, let us consider the function of two variables $\phi(a_1,a_2)$ that in some neighborhood of the point $(0,0)$ has the following behavior:
$$\phi(a_1,a_2) = \begin{cases}
                    -a_1(a_1+a_2) + 3a_1^2 +2a_2^2 , & \mbox{if } x\ge 0 \\
                    -a_1(a_1+a_2) + 10a_1^2 +2a_2^2, & \mbox{if } x<0.
                  \end{cases}$$
The function above has the gradient everywhere on $\R^2$. However, it is not twice differentiable at the local minimum $(0,0)$. Nevertheless, Assumption~\ref{assum:A5A} holds for $\ba' = (0,0)$, $C(\ba')$ equal to the $2\times 2$ identity matrix, and any $\ba\in\R^2$.}

\revi{Some works study the problem of escaping saddle points in centralized non-convex optimization \cite{nonconvexGE}, \cite{Anandkumar2016EfficientAF}, \cite{pemantle1990}.
However, it is shown in \cite{Anandkumar2016EfficientAF} that escaping from saddle points of the forth order and higher is an NP-hard problem, while for escaping from saddle points of a lower order the assumption on existence of second or third order derivatives of the objective function is needed. This work assumes objective functions to have a Lipschitz continuous gradient and does not require them to be twice differentiable.}

\section{Communication-based Algorithm}\label{subsec:comba}
In this subsection we consider a system, where agents can communicate and exchange information only with their neighbors. Besides such communication, the gradient of the local utility function is available for each agent. The neighbor set and the information which each agent can share are defined by a directed graph $G(t)$. We deal with the time-varying topology of the communication graph $G(t)$ that is $S$-strongly connected (see Definition\s\ref{def:CG}).
Another assumption, under which we will analyze the learning algorithm presented further, concerns the following continuous property of the gradient $\nabla\phi$:

\begin{assumption}\label{assum:CoordLip}
The gradient $\nabla \phi$ is coordinate-wise Lipschitz continuous on $\R^N$, i.e. there exists such positive constant $L$ that $\left|\frac{ \partial\phi(\ba_1)}{\partial a_i} - \frac{ \partial\phi(\ba_2)}{\partial a_i}\right|\le L\|\ba_1-\ba_2\|$ for all $i\in[N]$ and $\ba_1, \ba_2 \in \R^N$.
\end{assumption}

We remind that the goal of the agents is to find a local solution of the problem \eqref{eq:opt} by learning to play a local maximum of the potential function in the game $\Gamma(N,\boldsymbol A,\{U_i\},\phi)$. To reach this goal the agents utilize the \emph{push-sum protocol} introduced in Subsection\s\ref{subsec:pertps}.
%
%
At every moment of time $t\in \mathbb Z^+$ each agent $i$ maintains, besides her current action $a_i$,  vector   variables $\hat{\bba}_i(t) = (\hat{\ra}_i^1,\ldots,\hat{\ra}_i^N)$, \revi{which corresponds to the local joint action estimation}, $\bx_i(t)$, $\boldsymbol w_i(t)\in \mathbb{R}^N$, as  well  as  a  scalar  variable $y_i(t): y_i(0)=1$. These quantities are updated according to the \emph{communication-based procedure} described by Algorithm\s\ref{eq:ps}.

\begin{algorithm}
\caption{Communication-based algorithm}\label{eq:ps}
\begin{algorithmic}[1]
\State Let $y_i(0)=1$, $\bx_i(0)$ and $\boldsymbol w_i(0)\in \mathbb{R}^N$, $i\in[N]$, be some initial auxiliary variables, $t=0$.
\smallskip

\State Let $\{\gamma(t)\}_t$ be a specified sequence of time steps, $\fb_i(\hat{\bba}_i(t+1)) = (0,\ldots, \frac{\partial U_i}{\partial a_i}(\hat{\bba}_i(t+1)),\ldots, 0)^T\in\R^N$, and $\{\xib_i(t)\in\R^N\}_t$ be the sequence of random vectors.
\smallskip

\State $\boldsymbol w_i(t+1)=\sum_{j\in N^{in}_i(t)}\frac{\bx_j(t)}{d_j(t)}$.
\smallskip

\State $y_i(t+1)=\sum_{j\in N^{in}_i(t)}\frac{y_j(t)}{d_j(t)}$.
\smallskip

\State $\hat{\bba}_i(t+1)=\frac{\boldsymbol w_i(t+1)}{y_i(t+1)}$.
\smallskip

\State $\bx_i(t+1)=\boldsymbol w_i(t+1)+\gamma(t+1)[\fb_i(\hat{\bba}_i(t+1))+\xib_i(t)]$.
\smallskip

\State $a_i(t+1) = \hat{\ra}_i^i(t+1)$.
\smallskip

\State  $t:=t+1$, \textbf{goto} 3-7.
\end{algorithmic}
\end{algorithm}

Further, we assume that $\{\xib_i(t)\}_t$ is the sequence of independently and identically distributed (i.i.d.) random vectors satisfying the following assumption:

\begin{assumption}\label{assum:A4}
  The coordinates $\xi^k(t)$, $\xi^l(t)$ are independent, $\|\xi^k(t)\|$ is bounded almost surely and $\E(\xi^k(t))=0$ for all $t\in \mathbb Z^+$ and $k,l\in[N]$.
\end{assumption}

The process in Algorithm\s\ref{eq:ps} is a special case of the \emph{perturbed push-sum} algorithm, whose background and properties are discussed in Section\s\ref{subsec:pertps}.
In Algorithm\s\ref{eq:ps} we adapt the push-sum protocol \eqref{eq:gps} to the distributed learning of a potential function maximizer by letting $\ze_i(t+1)$ be $\gamma(t+1)[\fb_i(\hat{\bba}_i(t+1)) + \xib_i(t)]$.
The idea of adding a noise term $\xib_i(t)$ to the local optimization step 6 in Algorithm\s\ref{eq:ps} is similar
to one of simulated annealing \cite{bertsimas1993}. According to Assumption\s\ref{assum:A5A}, non-local optima are unstable with respect to the gradient descent direction of the global objective function and, thus, can be eliminated by the additional noise at the corresponding optimization step. Without this noise, the gradient procedure can guarantee convergence only to some critical point.
Under Assumptions\s\ref{assum:Ufgradient},\s\ref{assum:A4} and given that $\lim_{t\to\infty}\gamma(t)=0$, we conclude that $\lim_{t\to\infty}\|\ze_i(t)\|_1=0$ a.s. for all $i$.
The intuition behind this scheme is as follows: according to Theorem~\ref{th:th2}, the fact that $\lim_{t\to\infty}\|\ze_i(t)\|_1=0$ a.s. for all $i$, and Remark~\ref{rem:psprot}, the nodes' variables $\hat{\bba}_i(t+1)$ will track the average state $\bar{\bx}(t) = \frac{1}{N}\sum_{j=1}^{N}\bx_j(t)$ in the long run of Algorithm\s\ref{eq:ps}.
In the following we will demonstrate that, given Assumptions\s\ref{assum:A0}-\ref{assum:A4} and an appropriate choice of $\gamma(t)$ in Algorithm\s\ref{eq:ps}, the average state $\bar{\bx}(t)$ converges to some local potential function maximum $\ba^*$ almost surely, namely $\Pr\{\lim_{t\to\infty}\bar{\bx}(t)=\ba^*\}=1$, irrespective of the initial states. Thus, taking into account the last equality in Algorithm~\ref{eq:ps}, we can conclude that the joint action $\ba(t)$ converges to some local maximum of the potential function as $t$ tends to infinity.

\subsection{Communication-Based Algorithm: Convergence Result}
This subsection formulates the main convergence result concerning the \emph{communication-based algorithm}.
Firstly, we notice that under the rules of Algorithm\s\ref{eq:ps} the average of the vector variables $\bx_i(t), i\in[N]$, namely $\bbx(t) = \frac{1}{N}\sum_{j=1}^N\bx_j(t)$, is updated as follows (for more details see \cite{A1}):

\begin{align}\label{eq:ps_vecform}
\bbx(t+1) = \bbx(t)&+\gamma(t+1)\left[\frac{\nabla\phi(\bbx(t))}{N}+\Rb(t,\bbx(t))\right]\cr
&+ \gamma(t+1)\xib(t),
\end{align}
with
\[\Rb(t,\bbx(t)) = \frac{1}{N}\sum_{i=1}^N[\Fb(\bba(t+1)) - \nabla\phi(\bbx(t))],\]
where
\begin{align*}
\Fb(\bba(t+1))&=\sum_{i=1}^N\fb_i(\hat{\bba}_i(t+1)) \cr
&= \left(\frac{\partial U_1(\hat{\bba}_1(t+1))}{\partial a_1},\ldots,\frac{\partial U_N(\hat{\bba}_N(t+1))}{\partial a_N}\right)
\end{align*}
and
$\xib(t) = \frac1N\sum_{i=1}^N\xib_i(t)$. Further we will analyze the convergence properties of the procedure \eqref{eq:ps_vecform} to be able to describe the behavior of the communication-based procedure in Algorithm\s\ref{eq:ps}. The main result of this subsection is formulated in the following theorem.

\begin{thm}\label{th:cba}
 Let $\Gamma(N, \boldsymbol A, \{U_i\},\phi)$ be a potential game with $A_i=\mathbb R$.
 Let the parameter $\gamma(t)$ be such that $\gamma(t)=O\left(\frac{1}{t^{\nu}}\right)$, where $\frac{1}{2}<\nu\le 1$.
 Then, given $S$-strongly connected sequence of communication graphs and under Assumptions\s\ref{assum:Ufgradient}-\ref{assum:A4}, the sequence $\{\bbx(t)\}$ defined in \eqref{eq:ps_vecform}
  converges either to a point from the set of local maxima or to the boundary of one of its connected components almost surely irrespective of the initial states.
 Moreover, the joint action $\ba(t) = (a_1(t),\ldots, a_N(t))$ defined in Algorithm\s\ref{eq:ps} converges either to a point from the set of local maxima or to the boundary of one of its connected components almost surely.
\end{thm}

\begin{rem}
\revi{Note that connected or strongly connected graphs in consensus protocols are standard assumptions \cite{gossipNedic} \cite{A1}, \cite{NedicO09}, \cite{TAC}. Connectivity of communication graphs is required to guarantee sufficient ``mixing'' of information exchanged by agents.
In this paper, a time-varying topology is considered. In this case, an $S$-strongly connected sequence of communication graphs guarantees such information ``mixing''.}
\end{rem}

\begin{rem}
 The theorem above states almost sure convergence of the agents' join action to a local minimum of the objective function. To estimate the convergence rate in this case, one can use the result on the push-sum algorithm for distributed optimization provided in \cite{TAC}.
 By adapting the analysis in Theorem\s6 in that paper to the case of the communication-based procedure in \eqref{eq:ps_vecform}, one can demonstrate that this procedure converges sublinearly.
 \revi{Moreover, the proof of Theorem\s\ref{th:cba} follows similar argument as of Theorem\s5 in \cite{TAC}. The main difference consists in the setting of the push-sum algorithm. To apply this algorithm to the optimization problem \eqref{eq:opt} we added the gradient term $\gamma(t+1)\frac{\partial U_i}{\partial a_i}(\hat{\bba}_i(t+1))$ only to the update of the $i$th agent's action (see step 6 and definition of the function $\fb_i$, $i\in[N]$). The estimations for the actions of other agents are updated according to the standard push-sum consensus protocol (see steps 3-5).}
\end{rem}

\begin{rem}
\revi{Furthermore, if Assumption\s\ref{assum:A5A} does not hold, we can define the set $\boldsymbol A'$ to be the set of $\ba'\in \boldsymbol A_0$ for which there exists a symmetric positive definite matrix $C(\ba')$ such that $(\nabla\phi(\ba),C(\ba')(\ba-\ba'))\le 0$ for any $\ba\in U(\ba')$, where $U(\ba')$ is some open neighborhood of $\ba'$.  Then $\bbx(t)$ and $\ba(t)$ converge almost surely to a point from $\boldsymbol A_0\setminus \boldsymbol A'$.}
\end{rem}

\subsection{Proof of the Main Result}

Under Assumption\s\ref{assum:CoordLip} the gradient function $\nabla \phi(\bx)$ is a coordinate-wise Lipschitz function. This allows us to formulate the following lemma which we will need to prove Theorem\s\ref{th:cba}.

\begin{lem}\label{lem:lem1}
Let $\{\gamma(t)\}$ be a non-increasing sequence such that $\sum_{t=0}^{\infty}\gamma^2(t)<\infty$. Then, given $S$-strongly connected sequence of communication graphs and under Assumptions\s\ref{assum:Ufgradient}, \ref{assum:CoordLip}, and\s\ref{assum:A4}, there exists $c(t)$ such that the following holds a.s. for the process \eqref{eq:ps_vecform}:
 \begin{align*}
 \|\Rb(t,\bar{\bx}(t))&\|= \cr
 &\frac{1}{N}\|\sum_{i=1}^N[\Fb(\bba(t+1))-\nabla\phi(\bar{\bx}(t))]\|\le c(t),
 \end{align*}
 and $\sum_{t=0}^{\infty} \gamma(t+1)c(t)<\infty.$
\end{lem}
\begin{IEEEproof}
Since the functions $\nabla\phi$ is a coordinate-wise Lipschitz one (Assumption\s\ref{assum:CoordLip}) and taking into account \eqref{eq:pg1}, we conclude that for any $i$
\begin{align*}
&\|\sum_{i=1}^N[\Fb(\bba(t+1))-\nabla\phi(\bar{\bx}(t))]\|\cr
&\le\sum_{i=1}^{N}\left|\frac{\partial U_i(\hat{\bba}_i(t+1))}{\partial a_i}-\frac{\partial \phi(\bx(t))}{\partial a_i}\right|\cr
&=\sum_{i=1}^{N}\left|\frac{\partial \phi(\hat{\bba}_i(t+1))}{\partial a_i}-\frac{\partial \phi(\bx(t))}{\partial a_i}\right|\cr
&\quad\le L\sum_{i=1}^{N}\|\hat{\bba}_i(t+1)-\bbx(t)\|.
\end{align*}
Let $c(t)=\frac{L}{N}\sum_{i=1}^{N}\|\hat{\bba}_i(t+1)-\bbx(t)\|$. Then,
$\|\Rb(t,\bbx(t))\|\le c(t)$. Taking into account Assumptions\s\ref{assum:Ufgradient} and\s\ref{assum:A4}, we get that almost surely
\begin{align*}
\sum_{t=1}^{\infty}\gamma(t)\|\ze_i(t)\|_1 &= \sum_{t=0}^{\infty}\gamma^2(t)\|\fb_i(\hat{\bba}_i(t+1))+\xib_i(t)\|\cr
&\le k\sum_{t=0}^{\infty}\gamma^2(t)<\infty,
\end{align*}
where $k$ is some positive constant.
Hence, according to Theorem\s\ref{th:th2}(b),
\[
\sum_{t=0}^{\infty}\gamma(t+1)c(t)=\frac{L}{N}\sum_{i=1}^{N}\sum_{t=0}^{\infty}\|\hat{\bba}_i(t+1)-\bbx(t)\|\gamma(t+1)<\infty.
\]
\end{IEEEproof}

Now we are ready to prove Theorem~\ref{th:cba}.
\begin{IEEEproof}[Proof of Theorem\s\ref{th:cba}]
 To prove this statement we will use the general result formulated in Theorem\s\ref{th:th1} and Theorem\s\ref{th:th5}. The process \eqref{eq:ps_vecform} under consideration represents a special case of the recursive procedure \eqref{eq:basic}, where $\zbx(t) = \bbx(t)$, $\alpha(t)=\beta(t)=\gamma(t)$, $\fb(\zbx(t)) = \nabla\phi(\bbx(t))$, $\zq(t,\zbx(t)) = \Rb(t,\bbx(t))$, and $\Wb(t, \zbx(t),\omega) = \xib(t)$.
 In the following any $k_j$ denotes some positive constant.
 Note that, according to the choice of $\gamma(t)$, Assumption\s\ref{assum:Ufgradient}, and Lemma\s\ref{lem:lem1}, the conditions\s\eqref{th1:cond3} and \eqref{th1:cond4} of Theorem\s\ref{th:th1} hold. Thus, it suffices to show that there exists a sample function $V(t, \bx)$ of the process \eqref{eq:ps_vecform} satisfying the conditions \eqref{th1:cond1} and \eqref{th1:cond2} of Theorem\s\ref{th:th1}. A natural choice for such a function is the (time-invariant) function $V(t,\bx)=V(\bx)=-\phi(\bx)+C$, where $C$ is a constant chosen to guarantee the positiveness of $V(\bx)$ over the whole $\mathbb{R}^N$. Note that such constant always exists because of the continuity of $\phi$ and Assumption\s\ref{assum:Ufinfin}.
Then, $V(\bx)$ is non-negative and $V(\bx)\to\infty \quad \mbox{as} \quad \|\bx\|\to\infty$ (see again Assumption\s\ref{assum:Ufinfin}).
Now we show that the function $V(\bx)$ satisfies the condition \eqref{th1:cond2} in Theorem\s\ref{th:th1}.
Using the notation ${\tilde{\fb}}(t,\bbx(t)) = \frac{\nabla\phi(\bbx(t))}{N}+\Rb(t, \bbx(t))+\xib(t)$ and by the Mean-value Theorem applied to the function $V(\bx)$:
\begin{align*}
LV(\bx)& = \E V(\bbx(t+1)|\bbx(t)=\bx) - V(\bx) \cr
&= \gamma(t+1)\E (\nabla V(\tilde{\bx}),\tilde{\fb}(t,\bbx(t)))\cr
& = \gamma(t+1)\E(\nabla V(\bbx(t)),\tilde{\fb}(t,\bbx(t)))\cr
&\quad +\gamma(t+1)\E(\nabla V(\tilde{\bx})-\nabla V(\bbx(t)),\tilde{\fb}(t,\bbx(t))),
\end{align*}
where $\tilde{\bx}=\bar{\bx}(t)-\theta \gamma(t+1)\tilde{\fb}(t,\bbx(t))$ for some $\theta\in (0,1)$. Taking into account Assumption\s\ref{assum:CoordLip}, we
obtain that
\begin{align*}
 \|\nabla V(\tilde{\bx})-\nabla V(\bbx(t))\|&\le k_1 \gamma(t+1)\|\tilde{\fb}(t,\bbx(t))\|.
\end{align*}

Since
\begin{align*}
 \E(\nabla V(\bbx(t)),\tilde{\fb}(t,\bbx(t))) =& (\nabla V(\bbx(t)),\nabla\phi(\bbx(t)))\cr
 +&(\nabla V(\bbx(t)),\Rb(t,\bbx(t))),\cr
 \E\|\tilde{\fb}(t,\bbx(t))\|^2 =& \frac{1}{N^2}\|\nabla\phi(\bbx(t))\|^2 + \|\Rb(t,\bbx(t))\|^2\cr
 + \E\|\xib(t)\|^2 &+ \frac2N(\nabla\phi(\bbx(t)), \Rb(t,\bbx(t))),
\end{align*}
and due to the Cauchy-Schwarz inequality and Assumptions\s\ref{assum:Ufgradient},\s\ref{assum:A4}, we conclude that
\begin{align*}
LV(\bx)
& \le \frac{\gamma(t+1)}{N}(\nabla V(\bbx(t)),\nabla\phi(\bbx(t)))\cr
&\quad +k_2 \gamma(t+1)\|\Rb(t, \bbx(t))\| + k_3\gamma^2(t+1).
\end{align*}

Now recall that $\nabla V = -\nabla \phi$ and the function $V(\bx)$ is non-negative. Thus, we finally obtain that
\begin{align*}
L&V(\bx)\cr
&\le -\frac{\gamma(t+1)}{N}\|\nabla \phi(\bbx(t))\|^2+ g(t)(1+V(\bbx(t))),
\end{align*}
where
\begin{align*}
g(t) = k_2 \gamma(t+1)\|\Rb(t, \bbx(t))\| + k_3\gamma^2(t+1).
\end{align*}
Lemma\s\ref{lem:lem1} and the choice of $\gamma(t)$ imply that $g(t)>0$ and
$\sum_{t=0}^{\infty}g(t)<\infty$.
Thus, all conditions of Theorem\s\ref{th:th1} hold and, taking into account $\|\nabla \phi(\bbx(t))\|\in\Psi(\boldsymbol A_0)$, we conclude that either $\lim_{t\to\infty} \bbx(t)=\ba^*\in \boldsymbol A_0$ almost surely or $\bbx(t)$ converges to the boundary of a connected component of the set $\boldsymbol A_0$ almost surely.

Next we show that the process converges almost surely either to $\ba^*\in\boldsymbol A^*$ or to the boundary of a connected component of the set $\boldsymbol A^*$. To do this we use the result formulated in Theorem\s\ref{th:th5}. Obviously, under Assumptions\s\ref{assum:CoordLip} and\s\ref{assum:A4} conditions \eqref{th5:cond1} and \eqref{th5:cond2} in Theorem\s\ref{th:th5} hold. Thus, it suffices to check the finiteness of the sums in conditions (3) and (4) in Theorem\s\ref{th:th5}. Since $\gamma(t) = O(1/t^{\nu})$, $\nu\in(0.5,1]$, and
\[\sum_{k=t+1}^{\infty} \gamma^2(k)  = O(1/t^{2\nu-1})\quad \mbox{ (see Remark\s\ref{rem:series})},\]
we get
\[\sum_{t=0}^{\infty}\left(\frac{\gamma(t)}{\sqrt{\sum_{k=t+1}^{\infty}\gamma^2(k)}}\right)^3 = \sum_{t=0}^{\infty}O(1/t^{1.5})<\infty.\]

Moreover,
\begin{align*}
\|\Rb(t,\bbx(t))\|&\le c(t) = \frac{L}{N}\sum_{i=1}^{N}\|\hat{\bba}_i(t+1)-\bbx(t)\|
\end{align*}
 almost surely. Hence, according to Theorem\s\ref{th:th2}(a) and Assumptions\s\ref{assum:Ufgradient} and\s\ref{assum:A4},
$$ \|\Rb(t,\bbx(t))\|\le M \left(\lambda^t + \sum_{s=1}^{t}\lambda^{t-s}\gamma(s)\right)$$
almost surely for some positive constant $M$. Here the parameter $\lambda\in(0,1)$ is one defined in Theorem\s\ref{th:th2}.
Taking into account this and the fact that $\frac{\gamma(t)}{\sqrt{\sum_{k=t+1}^{\infty}\gamma^2(k)}}=O(\frac{1}{\sqrt{t}})$, we have
\begin{align*}
 \sum_{t=0}^{\infty}&\frac{\gamma(t)\|\Rb(t,\bbx(t))\|}{\sqrt{\sum_{k=t+1}^{\infty}\gamma^2(k)}}<\infty.
\end{align*}
The last inequality is obtained due to the fact that
$\sum_{t=0}^{\infty}\lambda^tO(\frac{1}{\sqrt t})<\infty$ and
\begin{align*}
\sum_{t=0}^{\infty}O\left(\frac{1}{\sqrt t}\right)&\sum_{s=1}^{t}\lambda^{t-s}\gamma(s)\cr
&\le\sum_{t=0}^{\infty}\sum_{s=1}^{t}\lambda^{t-s}O\left(\frac{1}{\sqrt s}\gamma(s)\right)<\infty,
\end{align*}
since $\sum_{s=1}^{\infty}\frac{1}{\sqrt s}\gamma(s)<\infty$ and, hence, one can use the result in \cite{RamVN09}, according to which any series of the type $\sum_{k=0}^{\infty}\left(\sum_{l=0}^{k}\beta^{k-l}\gamma_l\right)$ converges, if $\gamma_k\ge 0$ for all $k\ge0$, $\sum_{k=0}^{\infty}\gamma_k<\infty$, and $\beta\in(0,1)$.
%
Thus, all conditions of Theorem\s\ref{th:th5} are fulfilled.

It implies that the process $\bbx(t)$ converges almost surely either to $\ba^*\in\boldsymbol A^*$ or to the boundary of a connected component of the set $\boldsymbol A^*$.
On the other hand, by Theorem\s\ref{th:th2}(a), $\lim_{t\to\infty}\|\hat{\bba}_i(t)-\bbx(t)\|=0$ for all $i\in [N]$. It means that $\hat{\bba}_i(t)$ converges almost surely either to $\ba^*\in\boldsymbol A^*$ or to the boundary of a connected component of the set $\boldsymbol A^*$. Finally, due to the learning procedure in Algorithm\s\ref{eq:ps},
\[a_i(t+1) = \hat{\ra}_i^i(t+1).\]
Hence, irrespective of the initial states, the joint action $\ba(t)=(a_1(t),\ldots, a_N(t))$ converges almost surely either to $\ba^*\in\boldsymbol A^*$ or to the boundary of a connected component of the set $\boldsymbol A^*$.
\end{IEEEproof}

\section{Payoff-Based Algorithm}\label{sec:pba}
In this section we introduce the \emph{payoff-based algorithm} to learn a potential function maximizer in a potential game.
In contrast to the previous subsection, we deal here with a system, where at every moment of time $t$ each agent $i$ has only the access to her currently played action $a_i(t)$ and the experienced payoff $U_i(t)$, where $U_i(t) = U_i(a_1(t),\ldots,a_N(t))$. We introduce here the algorithm updating the agents' mixed strategies, which in their turn are chosen in a way to put more weight on the actions in the direction  of the utilities' ascent. We refer here to the idea of CALA (continuous action-set learning automaton) presented in the literature on learning automata \cite{Beigy}, \cite{Thatha} and use the Gaussian distribution as the basis for mixing agents' strategies. \revi{Note that Gaussian samplings are used in many works devoted to parameter estimations in unknown systems \cite{Thompson}, \cite{RePEc:cor:louvco:2011001}, \cite{LanNONConvex}. Random samplings allow a zero-order procedure to explore the system and to get an appropriate estimation of gradient's values.}

We consider a potential game $\Gamma(N,\boldsymbol A,\{U_i\},\phi)$, where each agent $i$ chooses her action from the action set $A_i = \R$ according to the normal distribution $\EuScript N(\mu_i,\sigma^2)$ with the density:
\begin{align}\label{eq:eq17}
 p_{\mu_i}(x)=\frac{1}{\sqrt{2\pi}\sigma}\exp\left\{-\frac{(x-\mu_i)^2}{2\sigma^2}\right\}.
\end{align}

The learning procedure updates the parameters $\mu_i$ for every $i\in[N]$ and $\sigma$ according to Algorithm\s\ref{eq:pba}.

\begin{algorithm}
\caption{Payoff-based algorithm}\label{eq:pba}
\begin{algorithmic}[1]
\State Let $\{\sigma(t)\}_t$ be a specified variance sequence, $\mu_i(0)$ be an initial mean parameter for the normal distribution of the agent $i$'s action, $a_i(0)\sim \EuScript N(\mu_i(0),\sigma^2(0))$, $U_i(0) = U_i(a_1(0),\ldots,a_N(0))$ be the value of the agent $i$'s utility function, $i\in[N]$, $t=0$.
\smallskip

\State Let $\{\gamma(t)\}_t$ be a specified sequence of time steps.
\smallskip

\State $\mu_i(t+1) = \mu_i(t) +\gamma(t+1)\sigma^3(t+1)\left[U_i(t)\frac{a_i(t) -\mu_i(t)}{\sigma^2(t)} \right].$
\smallskip

\State $a_i(t+1) \sim \EuScript N(\mu_i(t+1),\sigma^2(t+1))$.
\smallskip

\State $t:=t+1$, \textbf{goto} 3-4.
\end{algorithmic}
\end{algorithm}

Thus, the payoff-based algorithm is of stochastic nature itself. The agents choose their actions according to the probabilistic distributions. The parameters of these distributions are updated iteratively according to the actually obtained payoff information.


In the case of the payoff-based algorithm above, we will need a less restrictive assumption on the continuity of the gradient $\nabla\phi$.

\begin{assumption}\label{assum:Lip}
The gradient $\nabla \phi$ is Lipschitz continuous on $\R^N$, i.e. there exists such a positive constant $L$ that $\|\nabla\phi(\ba_1) - \nabla\phi(\ba_2)\|\le L\|\ba_1-\ba_2\|$ for any $\ba_1, \ba_2 \in \R^N$.
\end{assumption}

However, for the analysis of the payoff-based algorithm we make the following assumption regarding the behavior of utility functions at infinity.

\begin{assumption}\label{assum:infbeh}
The utility functions $U_i(\ba)$, $i\in[N]$, are continuous on the whole $\R^N$ and their absolute values grow not faster than a linear function of $\ba$, given $\ba$ has a large norm, i.e. there exist such constants $K>0$ and $c>0$ that
$|U_i(\ba)|\le c(1 + \|\ba\|)$ for any $i\in[N]$, if $\|\ba\|\ge K$.
\end{assumption}

In the following subsection we will provide the analysis of Algorithm\s\ref{eq:pba}.
\revi{To the best of the author's knowledge, no payoff-based algorithm has been presented in the literature so far that would be applicable to distributed non-concave optimization \eqref{eq:opt} in multi-agent systems. For example, the known technique on a smooth approximation of objective functions introduced by Nesterov in \cite{RePEc:cor:louvco:2011001} suffers from a shortcoming that two points are needed to estimate the gradient of the smoothed function. Thus, agents would need an \emph{extra coordination} to get an appropriate estimation, if this technique is applied to multi-agent optimization. Moreover, Nesterov's approach assumes a so called smoothing parameter to be fixed and predefined. This parameter defines how precisely the smoothed function approximates the objective one \cite{LanNONConvex}.
As the result, algorithms based on this approach may converge only to a \emph{neighborhood} of a critical point \cite{LanNONConvex}.
To rectify these issues, this work uses $\frac{U_i(a_i-\mu_i)}{\sigma^2}$,
$a_i\sim\EuScript N(\mu_i,\sigma^2)$, $U_i = U_i(\ba)$ to obtain an estimation of gradients. This estimation is based only on local utilities' values and, hence, no extra coordination between agents is required.}
In the following subsection we will show that under Assumptions\s\ref{assum:A0}-\ref{assum:Ufinfin},\s\ref{assum:Lip},\s\ref{assum:infbeh} and an appropriate choice of the variance parameter $\sigma(t)$ and the step-size parameter $\gamma(t)$, the joint actions updated according to Algorithm\s\ref{eq:pba} converge in probability to a local maximum of the potential function as time tends to infinity.
\vspace{-0.3cm}

\subsection{Payoff-Based Algorithm: Convergence Result}
In this subsection we consider the \emph{payoff-based procedure} in Algorithm\s\ref{eq:pba} and formulate the main result concerning its convergence.
Firstly, we calculate the expectation of the term $\frac{U_i(t)(a_i(t) - \mu_i(t))}{\sigma^2(t)}$.
Let us consider the following notations: $\bx=(x_1,\ldots,x_N)$, $p_{\bmu}(\bx) = \prod_{i\in [N]} p_{\mu_i}(x_i)$, and $p_{\bmu_{-i}}(\bx_{-i}) = \prod_{j\ne i} p_{\mu_j}(x_j)$, where $p_{\mu_i}$, $i\in[N]$, is introduced in \eqref{eq:eq17}.
Taking into account Assumptions~\ref{assum:Ufgradient},\s\ref{assum:Ufinfin},~\ref{assum:infbeh}, and the fact that the coordinates of the joint action $\ba(t)$ are independently distributed and the $i$th coordinate has the distribution $\EuScript N_i$=$\EuScript N(\mu_i(t), \sigma^2(t))$, we obtain\footnote{For the sake of notation simplicity, we omit the argument
$t$ in some estimation.}:
\begin{align}\label{eq:expectation}
 \E&\{\frac{U_i(a_i - \mu_i)}{\sigma^2}|a_i\sim\EuScript N(\mu_i,\sigma), i\in[N]\}  \cr
 &= \frac{1}{\sigma^2}\int_{\R^N}U_i(\bx)(x_i - \mu_i)p_{\bmu}(\bx)d\bx\cr
 &= -\int_{\R^N}U_i(\bx)  p_{\bmu_{-i}}(\bx_{-i}) d\left(e^{-\frac{(x_i-\mu_i)^2}{2\sigma^2}}\right) d\bx_{-i}\cr
 &=-\int_{\R^{N-1}}\left(U_i(\bx)e^{-\frac{(x_i-\mu_i)^2}{2\sigma^2}}\right)\bigg|_{x_i =-\infty}^{\infty}p_{\bmu_{-i}}(\bx_{-i})d\bx_{-i}\cr
 &\quad+\int_{\R^N}\frac{\partial U_i(\bx)}{\partial x_i}p_{\bmu}(\bx)d\bx\cr
 &=\int_{\R^N}\frac{\partial U_i(\bx)}{\partial x_i}p_{\bmu}(\bx)d\bx.
\end{align}
For each $\bmu\in\R^N$ let us introduce the vector $\nabla\tilde{\phi}(\bmu)$ with the $i$th coordinate equal to $\int_{\R^N}\frac{\partial \phi(\bx)}{\partial x_i}p_{\bmu}(\bx)d\bx$. According to equality \eqref{eq:pg1}, we conclude that
$$\int_{\R^N}\frac{\partial \phi(\bx)}{\partial x_i}p_{\bmu}(\bx)d\bx = \int_{\R^N}\frac{\partial U_i(\bx)}{\partial x_i}p_{\bmu}(\bx)d\bx.$$
Taking this equality and equality \eqref{eq:expectation} into account, we can rewrite Algorithm\s\ref{eq:pba} as follows:
\begin{align}\label{eq:pba_vecform}
 \bmu(t&+1)=\cr
 &\bmu(t)+ \sigma^3(t+1)\gamma(t+1)[\nabla{\phi}(\bmu(t)) +\Qb(\bmu(t))]\cr
 &\qquad \quad+ \gamma(t+1)[\sigma^3(t+1)\Mb(\ba,t,\bmu(t))],
\end{align}
where $\Qb(\bmu(t))=\nabla\tilde{\phi}(\bmu(t)) - \nabla{\phi}(\bmu(t))$ and $\Mb(\ba,t,\bmu(t))$ is the $N$-dimensional vector with the coordinates
\begin{align*}
M_i(\ba,t,\bmu(t)) = &U_i(t)\frac{a_i(t) -\mu_i(t)}{\sigma^2(t)} \cr
&- \int_{\R^N}\frac{\partial U_i(\bx)}{\partial x_i}p_{\bmu}(\bx)d\bx, \quad i\in[N].
\end{align*}
Due to \eqref{eq:expectation},
\begin{align}\label{eq:martdif}
 \E \{\Mb(\ba,t, \bmu(t))|a_i\sim\EuScript N(\mu_i(t),\sigma(t)), i\in[N]\} =\boldsymbol 0.
\end{align}

Now we are ready to formulate the main result on the payoff-based procedure in Algorithm\s\ref{eq:pba} and its vector form \eqref{eq:pba_vecform}.

\begin{thm}\label{th:mthpba}
 Let $\Gamma(N, \boldsymbol A, \{U_i\},\phi)$ be a potential game with $A_i=\mathbb R$.
 Let the parameter $\gamma(t)$ and $\sigma(t)$ be such that $\gamma(t)>0$, $\sigma(t)>0$,
 \begin{enumerate}[(1)]
	\item\label{mthpba:cond1} $\sum_{t=0}^{\infty}\gamma(t)\sigma^3(t) = \infty$, $\sum_{t=0}^{\infty}\gamma(t)\sigma^4(t) < \infty$,
	\medskip
	
	\item\label{mthpba:cond2}  $\sum_{t=0}^{\infty}\gamma^2(t)<\infty$, $\sum_{t=0}^{\infty}\left(\frac{\gamma(t)}{\sqrt{\sum_{k=t+1}^{\infty}\gamma^2(k)}}\right)^3<\infty$,
	\medskip
	
	\item\label{mthpba:cond4}  $\sum_{t=0}^{\infty}\frac{\gamma(t)\sigma^4(t)}{\sqrt{\sum_{k=t+1}^{\infty}\gamma^2(k)}}<\infty$.
	\medskip
 \end{enumerate}

 Then under Assumptions\s\ref{assum:Ufgradient}-\ref{assum:Ufinfin} and\s\ref{assum:Lip},\ref{assum:infbeh} the sequence $\{\bmu(t)\}$ defined by Algorithm\s\ref{eq:pba} converges either to a point from the set of local maxima or to the boundary of one of its connected components almost surely irrespective of the initial parameters.
 Moreover, the agents' joint action updated by the payoff-based procedure defined in Algorithm\s\ref{eq:pba} converges in probability either to a point from the set of local maxima or to the boundary of one of its connected components.
\end{thm}

\begin{rem}\label{rem:series}
 There exist such sequences $\{\gamma(t)\}$, $\{\sigma(t)\}$ that the conditions \eqref{mthpba:cond1}-\eqref{mthpba:cond4} of Theorem\s\ref{th:mthpba} hold. For example, let us consider $\gamma(t) = 1/t^{0.6}$ and $\sigma(t) = 1/t^{0.13}$. Notice that in this case
 \begin{align*}
 \sum_{k=t+1}^{\infty} \gamma^2(k) = \sum_{k=t+1}^{\infty} \frac{1}{k^{1.2}}\sim\int_{t+1}^{\infty}\frac{1}{x^{1.2}}dx = O(1/t^{0.2}).
\end{align*}
Hence,
\[\sum_{t=0}^{\infty}\left(\frac{\gamma(t)}{\sqrt{\sum_{k=t+1}^{\infty}\gamma^2(k)}}\right)^3 = \sum_{t=0}^{\infty}O(1/t^{1.5})<\infty,\]

\[\sum_{t=0}^{\infty}\frac{\gamma(t)\sigma^4(t)}{\sqrt{\sum_{k=t+1}^{\infty}\gamma^2(k)}} = \sum_{t=0}^{\infty}O(1/t^{1.02})<\infty.\]
\end{rem}

\begin{rem}
 To estimate the convergence rate of the proposed payoff-based procedure, one can apply the standard technique proposed in (5.292), \cite{stochprogr}, to estimate $\|\bmu(t)-\bmu^*\|$ based on \eqref{eq:pba_vecform}. The conjecture is that the convergence rate is sublinear, at least in the case of strongly concave objective functions (see also Theorem\s10 in \cite{TatKamg}).
\end{rem}

\begin{rem}
 \revi{Furthermore, if Assumption\s\ref{assum:A5A} does not hold, we can define the set $\boldsymbol A'$ to be the set of $\ba'\in \boldsymbol A_0$ for which there exists a symmetric positive definite matrix $C(\ba')$ such that $(\nabla\phi(\ba),C(\ba')(\ba-\ba'))\le 0$ for any $\ba\in U(\ba')$, where $U(\ba')$ is some open neighborhood of $\ba'$. Then $\bmu(t)$ and $\ba(t)$ converge to a point from $\boldsymbol A_0\setminus \boldsymbol A'$ almost surely and in probability respectively.}
\end{rem}

\begin{IEEEproof}
We will deal with the equivalent formulation \eqref{eq:pba_vecform} of the learning procedure in Algorithm\s\ref{eq:pba}.
In the following any $k_j$ denotes some positive constant.
To prove the claim we will demonstrate that the process \eqref{eq:pba_vecform} fulfills the conditions in Theorem~\ref{th:th1} and\s\ref{th:th5}, where $d=N$, $\zbx(t)=\bmu(t)$, $\alpha(t)=\gamma(t)\sigma^3(t)$, $\beta(t)=\gamma(t)$, $\fb(\zbx(t)) = \nabla\phi(\bmu)$, $\zq(t,\zbx(t))=\Qb(\bmu(t))$, and $\Wb(t,\zbx(t),\omega)=\sigma^3(t)\Mb(\ba,t,\bmu(t))$.
First of all, we show that there exists a sample function $V(t, \bmu)$ of the process \eqref{eq:pba_vecform} satisfying conditions\s\eqref{th1:cond1} and\s\eqref{th1:cond2} of Theorem~\ref{th:th1}. Let us consider the following time-invariant function
$$V(\bmu) = -\phi(\bmu)  + \sum_{i=1}^N h(\mu_i) + C ,$$
where\footnote{The constant $K$ below is one from Remark~\ref{rem:linbeh}.}
$$h(\mu_i) = \begin{cases}
              (\mu_i - K)^2, &\mbox{ if } \mu_i\ge K;\\
               0, &\mbox{ if } |\mu_i|\le K;\\
              (\mu_i + K)^2, &\mbox{ if } \mu_i\le -K;\\
             \end{cases}
$$
and $C$ is chosen in such way that $V(\bmu)>0$ for all $\bmu\in\R^N$. Thus, $V(\bmu)$ is po\-sitive on $\R^N$ and $\lim_{\|\bmu\|\to\infty} V(\bmu) = \infty$.

Further we will use the notation $\E_{\bmu,\sigma(t)}\{\xi(\ba)\}$ for the expectation of any random variable $\xi$ dependent on $\ba$ given that $\ba$ has a normal distribution with the parameters $\bmu,\sigma^2$, namely $\E_{\bmu,\sigma(t)}\{\xi(\ba)\} = \E\{\xi(\ba)|\ba\sim\EuScript N(\bmu,\sigma^2)\}$.

 We can use the Mean-value Theorem for the function $V(\bmu)$ to get
\begin{align}\label{eq:LV}
LV(\bmu)& = \E_{\bmu,\sigma(t)} V({\bmu} + \sigma^3(t+1)\gamma(t+1)\tilde{\fb}(\ba,t,{\bmu}))\cr
&\qquad\qquad\qquad\qquad\qquad\qquad - V({\bmu})\cr
& = \E_{\bmu,\sigma(t)}\sigma^3(t+1)\gamma(t+1)(\nabla V(\tilde{\bmu}),\tilde{\fb}(\ba,t,{\bmu}))\cr
& =\sigma^3(t+1)\gamma(t+1)\E_{\bmu,\sigma(t)} \{(\nabla V({\bmu}),\tilde{\fb}(\ba,t,{\bmu}))\cr
&\quad+(\nabla V(\tilde{\bmu})-\nabla V({\bmu}),\tilde{\fb}(\ba,t,{\bmu}))\},
\end{align}
where
\begin{align*}
  \tilde{\fb}(\ba,t,{\bmu}) &= \nabla{\phi}(\bmu) +\Qb(\bmu) + \Mb(\ba,t,\bmu),\cr
  \tilde{\bmu} &= {\bmu} + \theta\sigma^3(t+1)\gamma(t+1)\tilde{\fb}(\ba,t,{\bmu})
\end{align*}
for some $\theta\in(0,1)$. We proceed by estimating the terms in \eqref{eq:LV}. Let $h(\bmu)=\sum_{i=1}^N h_i(\bmu)$. Then, taking into account \eqref{eq:martdif} and the fact that the vector $\nabla h(\bmu)$ has coordinates that are linear in $\bmu$, we get:
\begin{align*}
 \E_{\bmu,\sigma(t)}& \{(\nabla V({\bmu}),\tilde{\fb}(\ba,t,{\bmu}))\} \cr
 &= -(\|\nabla{\phi}(\bmu)\|^2 - (\nabla{\phi}(\bmu),\nabla{h}(\bmu)) ) \cr
 &\quad+(\nabla h(\bmu) - \nabla{\phi}(\bmu), \Qb(\bmu))\cr
 &\le-(\|\nabla{\phi}(\bmu)\|^2- (\nabla{\phi}(\bmu),\nabla{h}(\bmu)) ) \cr
 &\quad+ k_1\|\Qb(\bmu)\|(1+V(\bmu)),
 \end{align*}
where  the last inequality is due to the Cauchy-Schwarz inequality and Assumption~\ref{assum:Ufgradient}. Thus, using again the Cauchy-Schwarz inequality, Assumption~\ref{assum:Lip}, and the fact that
$\nabla V(\bmu)$ is Lipschitz continuous, we obtain from \eqref{eq:LV}:
\begin{align}\label{eq:LV1}
&LV(\bmu,t)\le-\sigma^3(t+1)\gamma(t+1)\cr
&\quad\times(\|\nabla{\phi}(\bmu)\|^2 -(\nabla{\phi}(\bmu),\nabla{h}(\bmu)))\cr
&\qquad\qquad+\sigma^3(t+1)\gamma(t+1) k_2(1+V(\bmu))\|\Qb(\bmu)\|\cr
&\quad + \sigma^6(t+1)\gamma^2(t+1)\E_{\bmu,\sigma(t)}\|\tilde{\fb}(\ba, t,{\bmu}))\|^2.
\end{align}
Recall the definition of the vector $\nabla \tilde{\phi}(\bmu)$:
\begin{align}\label{eq:gradmix}
 \nabla \tilde{\phi}(\bmu) = \int_{\R^N}\nabla \phi(\bx)p_{\bmu}(\bx)d\bx.
\end{align}
Since $\Qb(\bmu(t))=\nabla\tilde{\phi}(\bmu(t)) - \nabla{\phi}(\bmu(t))$ and due to Assumptions~\ref{assum:Ufgradient},\s\ref{assum:Lip} and \eqref{eq:gradmix}, we can write the following:
\begin{align}\label{eq:Qterm}
 \|&\Qb(\bmu(t))\| \le \int_{\R^N}\|\nabla\phi(\bmu) - \nabla \phi(\bx)\| p_{\bmu}(\bx) d\bx \cr
 &\le \int_{\R^N} L \|\bmu - \bx\| p_{\bmu}(\bx) d\bx \cr
 &\le \int_{\R^N} L \left(\sum_{i=1}^N|\mu_i - x^i|\right) p_{\bmu}(\bx) d\bx= O(\sigma(t)),
\end{align}
Due to Assumption~\ref{assum:infbeh} and \eqref{eq:expectation},
\begin{align}\label{eq:variance}
 \sigma^6(t+1)\E_{\bmu,\sigma(t)}\|\Mb(\ba,t,\bmu)\|^2 = f_1(\bmu),
\end{align}
where $f_1$ depends on $\bmu$ and is bounded by a quadratic function. Hence, we conclude that
$$\sigma^6(t+1)\E_{\bmu,\sigma(t)}\|\tilde{\fb}(\ba,t,{\bmu}))\|^2\le k_3(1+V(\bmu)).$$
Thus, \eqref{eq:LV1} implies:
\begin{align}\label{eq:finalineq}
 &LV(\bmu,t) \cr
 &\le -\sigma^3(t+1)\gamma(t+1)(\|\nabla{\phi}(\bmu)\|^2 -(\nabla{\phi}(\bmu),\nabla{h}(\bmu)))\cr
 &\qquad\qquad\qquad\qquad\qquad\qquad\quad+g(t)(1+V(\bmu)),
\end{align}
where $g(t)=O(\sigma^4(t)\gamma(t)+\gamma^2(t))$, i.e. $\sum_{t=1}^{\infty}g(t)<\infty$, according to the choice of the sequence $\gamma(t)$ and $\sigma(t)$ (see conditions\s\eqref{mthpba:cond1} and\s\eqref{mthpba:cond2} in the theorem formulation). Note also that, according to the definition of the function $h$, $\|\nabla{\phi}(\bmu)\|^2 -(\nabla{\phi}(\bmu),\nabla{h}(\bmu))\ge 0$, where equality holds only on critical points of the function $\phi$ (see Remark\s\ref{rem:linbeh}). Thus, conditions\s\eqref{th1:cond1} and\s\eqref{th1:cond2} of Theorem~\ref{th:th1} hold.
Conditions\s\eqref{th1:cond3} and\s\eqref{th1:cond4} of Theorem~\ref{th:th1} hold, due to \eqref{eq:Qterm} and Assumption\s\ref{assum:Ufgradient} respectively.
Moreover, taking into account Theorem~\ref{th:finiteness} and \eqref{eq:finalineq}, we conclude that the norm $\|\bmu(t)\|$ is bounded almost surely for all $t$. Hence, condition \eqref{th1:cond5} of Theorem~\ref{th:th1} holds as well.
Thus, all conditions of Theorem~\ref{th:th1} are fulfilled.
It implies that $\lim_{t\to\infty}\bmu(t) = \bmu^*$ almost surely, where $\bmu^*$ is some critical point of $\phi$, or $\bmu(t)$ converges to the boundary of a connected component of the set $\boldsymbol A_0$.
Moreover, since $\sigma(t)\to 0 $ as $t\to\infty$, we conclude that $\ba(t)$ converges to $\bmu^*$ or to the boundary of a connected component of the set $\boldsymbol A_0$ in probability irrespective of the initial states.

Further we verify the fulfillment of conditions in Theorem~\ref{th:th5} to prove that this critical point $\ba^*=\bmu^*$ is necessarily a local maximum of $\phi$. Let $\bmu'$ denote a critical point of the function $\phi$ that is not in the set of local maxima $\boldsymbol A^*$.
We show that there exists some $\delta'>0$ such that $|\Tr[A(t,\bmu)-A(t,\bmu')]|\le k_4\|\bmu - \bmu'\|$, for any $\bmu$: $\|\bmu - \bmu'\|<\delta'$, where $A_{ii}(t,\bmu)=\E_{\bmu,\sigma(t)}\sigma^6(t+1)M_i^2(\ba,t,\bmu)$. Indeed,
\begin{align}\label{eq:trace}
 &\frac{1}{\sigma^6(t+1)}|\Tr[A(t,\bmu)-A(t,\bmu')]|\cr
 &=\left|\sum_{i=1}^N \E_{\bmu,\sigma(t)} M_i^2(\ba,t,\bmu) - \E_{\bmu',\sigma} M_i^2(\ba,t,\bmu')\right|\cr
 & \le | \sum_{i=1}^N\E_{\bmu,\sigma(t)} U^2_i(t)\frac{(a_i(t) -\mu_i)^2}{\sigma^4(t)} \cr
 &\qquad\qquad\qquad\qquad- \E_{\bmu',\sigma} U^2_i(t)\frac{(a_i(t) -\mu'_i)^2}{\sigma^4(t)}|\cr
  & \quad \qquad\qquad\qquad\qquad + \|\nabla \tilde{\phi}(\bmu')\|^2 - \|\nabla \tilde{\phi}(\bmu)\|^2.
\end{align}
Since $\nabla \phi$ is Lipschitz continuous, we can use \eqref{eq:gradmix} to get
\begin{align*}
 &\|\nabla\tilde{\phi}(\bmu) - \nabla \tilde{\phi}(\bmu')\|\cr
 &= \|\int_{\R^N}(\nabla \phi(\bx)p_{\bmu}(\bx) - \nabla \phi(\bx)p_{\bmu'}(\bx))d\bx\|\cr
 &\le\int_{\R^N}\|\nabla\phi(\by + \bmu) - \nabla\phi(\by + \bmu')\|p(\by)d\by \cr
 &\le k_5\|\bmu-\bmu'\|,
 \end{align*}
 where $p(\by) = \frac{1}{(2\pi \sigma^2)^{N/2}} e^{-\frac{Ny^2}{2\sigma^2}}$.
Hence, according to Assumptions~\ref{assum:Ufgradient},\s\ref{assum:Lip} and the last inequality,
\begin{align}\label{eq:trace2}
 \|\nabla \tilde{\phi}(\bmu')\|^2 - \|\nabla \tilde{\phi}(\bmu)\|^2&\le k_6(\|\nabla \tilde{\phi}(\bmu')\| - \|\nabla \tilde{\phi}(\bmu)\|)\cr
 &\le k_6\|\nabla \tilde{\phi}(\bmu') - \nabla \tilde{\phi}(\bmu)\|\cr
 &\le k_7\|\bmu-\bmu'\|.
\end{align}
Moreover,
\begin{align}\label{eq:trace3}
 \sigma^6\E_{\bmu,\sigma(t)}[\{ U^2_i\frac{(a_i -\mu_i)^2}{\sigma^4}\}&-\sigma^6 \{ U^2_i\frac{(a_i -\mu'_i)^2}{\sigma^4}\}]\cr
  &= u_i(t,\bmu) - u_i(t,\bmu'),
\end{align}
where
\begin{align}\label{eq:difunderint}
u_i(\bmu) = \int_{\R^N}\sigma^2U^2_i(\bx){(x^i -\mu_i)^2}p_{\bmu}(\bx)d\bx.
\end{align}
The function above is Lipschitz continuous in some $\delta$-neighborhood of $\bmu'$, since there exists $\delta>0$ such that the gradient  $\nabla u_i(\bmu)$ is bounded for any $\bmu$: $\|\bmu - \bmu'\|<\delta$. Indeed, due to Assumption~\ref{assum:A0} and if $\|\bmu - \bmu'\|<\delta$, the mean vectors $\bmu'$ and $\bmu$ in \eqref{eq:trace3} are bounded. Next, taking into account the behavior of each $U_i(\bx)$, when $\bx\to\infty$ (see Assumption~\ref{assum:infbeh}), we can use finiteness of moments of a normal random vector with a
bounded mean vector and apply the sufficient condition for uniform convergence of integrals based on majorants (see \cite{zorich}, Chapter 17.2.3.) to conclude that the integral in \eqref{eq:difunderint} can be differentiated under the integral sign with respect to the parameter $\mu^j$ for any $j\in[N]$. For the same reason of moments' finiteness, the partial derivative $\frac{\partial u_i(\bmu)}{\partial \mu^j}$ is bounded for all $i,j\in[N]$. According to the Mean-value Theorem for each function $u_i$, $i\in[N]$, this implies that
\begin{align}\label{eq:trace4}
 \sum_{i=1}^N|u_i(\bmu)-u_i(\bmu')|\le k_8\|\bmu - \bmu'\|.
\end{align}
Substituting \eqref{eq:trace2} - \eqref{eq:trace4} into \eqref{eq:trace} and taking into account that, according to Assumption~\ref{assum:Lip}, $\|\nabla \phi(\bmu) \|^2 = \|\nabla \phi(\bmu) - \nabla \phi(\bmu')\|^2\le L^2\|\bmu-\bmu'\|^2$ for all $\bmu\in\R^N$, we obtain that there exists $\delta'\le\delta$ such that for any $\bmu:$ $\|\bmu - \bmu'\|<\delta'$
\begin{align*}
 \|\nabla \phi(\bmu)\|^2 + |\Tr[A(t,\bmu)-A(t,\bmu')]|\le k_8\|\bmu - \bmu'\|.
\end{align*}
Thus, condition (\ref{th5:cond1}) of Theorem~\ref{th:th5} holds.

Since $\|\bmu'\|<\infty$ (Assumption~\ref{assum:A0}) and due to Assumptions~\ref{assum:Ufgradient},\s\ref{assum:Ufinfin}, and Remark~\ref{rem:linbeh}, $\sigma^{12}\E_{\bmu,\sigma(t)}\|\Mb(\ba,t,\bmu)\|^4<\infty$ for any $\bmu:$ $\|\bmu - \bmu'\|<\delta'$. Hence, condition (\ref{th5:cond2}) of Theorem~\ref{th:th5} is also fulfilled.

Finally, taking into account the choice of the sequences $\{\gamma(t)\}$, $\{\sigma(t)\}$ (see conditions \eqref{mthpba:cond1}-\eqref{mthpba:cond4}) and the estimation \eqref{eq:Qterm}, we conclude that the last two conditions of Theorem~\ref{th:th5} are also fulfilled. It allows us to conclude that $\Pr\{\lim_{t\to\infty}\bmu(t)\in \boldsymbol A_0\setminus \boldsymbol A^*\}=0$.  As $\lim_{t\to 0}\sigma(t)=0$, the players' joint actions $\ba(t)$, chosen according to the rules in Algorithm\s\ref{eq:pba}, converge in probability to a local maximum of the potential function $\phi$ or to the boundary  of a connected component of the set $\boldsymbol A^*$ irrespective of the initial states.
\end{IEEEproof}

\subsection{Payoff-based Algorithm Based on Two-Point Evaluations}\label{subsec:pb2p}
\revi{In this subsection, we present a version of the payoff-based algorithm introduced above, which is based on evaluations of utility functions in two points. For this purpose let us replace each estimation of the gradients in mixed strategies (see \eqref{eq:expectation}), namely $U_i(t)\frac{a_i(t) -\mu_i(t)}{\sigma^2(t)}$, by the expression $(U_i(t) - U_i(\bmu(t)))\frac{a_i(t) -\mu_i(t)}{\sigma^2(t)}$. Thus, the modified procedure can be described by Algorithm 3.}
\begin{algorithm}\label{eq:pba2p}
\caption{Two-point payoff-based algorithm}
\begin{algorithmic}[1]
\State Let $\{\sigma(t)\}_t$ be a specified variance sequence, $\mu_i(0)$ be an initial mean parameter for the normal distribution of the agent $i$'s action, $a_i(0)\sim \EuScript N(\mu_i(0),\sigma^2(0))$, $U_i(0) = U_i(a_1(0),\ldots,a_N(0))$ be the value of the agent $i$'s utility function, $i\in[N]$, $t=0$.
\smallskip

\State Let $\{\gamma(t)\}_t$ be a specified sequence of time steps.
\smallskip

\State \begin{align*}
         \mu_i(t+1) &= \mu_i(t) +\gamma(t+1)\sigma^3(t+1)\cr
         &\times\left[(U_i(t) - U_i(\bmu(t)))\frac{a_i(t) -\mu_i(t)}{\sigma^2(t)} \right]
       \end{align*}

\State $a_i(t+1) \sim \EuScript N(\mu_i(t+1),\sigma^2(t+1))$.
\smallskip

\State $t:=t+1$, \textbf{goto} 3-4.
\end{algorithmic}
\end{algorithm}

\revi{Due to the fact that
\[\E\{U_i(\bmu(t))\frac{a_i(t) -\mu_i(t)}{\sigma^2(t)}| a_i(t)\sim\EuScript N(\mu_i(t),\sigma^2(t))\} = 0, \]
we conclude that, analogously to $U_i(t)\frac{a_i(t) -\mu_i(t)}{\sigma^2(t)}$ (see \eqref{eq:expectation}), $(U_i(t) - U_i(\bmu(t)))\frac{a_i(t) -\mu_i(t)}{\sigma^2(t)}$ is an unbiased estimation of the gradients in mixed strategies. However, this estimation has a bounded variance irrespective of the mean value $\bmu(t)$, in contrast to $U_i(t)\frac{a_i(t) -\mu_i(t)}{\sigma^2(t)}$ (see \eqref{eq:variance}). Indeed,
according to Assumption\s\ref{assum:infbeh} and by using Mean Value Theorem, we conclude that there exists a constant $k>0$
\begin{align}\label{eq:variance2p}
&\E_{\bmu,\sigma(t)}\frac{(U_i(\ba) - U_i(\bmu))^2(a_i -\mu_i)^2}{\sigma^4} \cr
&\quad= k\int_{\R^N}\frac{(\bx - \bmu)^2(x_i -\mu_i)^2}{\sigma^4(t)}p_{\bmu}(\bx)d\bx = Nk.
\end{align}
}
\revi{Note that due to the equality above, the proof of Theorem\s\ref{th:mthpba} can be repeated for Algorithm 3, where $g(t) = O(\sigma^4(t)\gamma(t)+\gamma^2(t))$ in \eqref{eq:finalineq} in  is replaced by $g(t) = O(\sigma^4(t)\gamma(t))$. Thus, we obtain the following theorem.}
\begin{thm}\label{th:mthpba2p}
 Let $\Gamma(N, \boldsymbol A, \{U_i\},\phi)$ be a potential game with $A_i=\mathbb R$.
 Let the parameter $\gamma(t)$ and $\sigma(t)$ be such that $\gamma(t)>0$, $\sigma(t)>0$,
 \begin{enumerate}[(1)]
	\item\label{mthpba:cond1} $\sum_{t=0}^{\infty}\gamma(t)\sigma^3(t) = \infty$, $\sum_{t=0}^{\infty}\gamma(t)\sigma^4(t) < \infty$,
	\medskip
	
	\item\label{mthpba:cond2} $\sum_{t=0}^{\infty}\left(\frac{\gamma(t)}{\sqrt{\sum_{k=t+1}^{\infty}\gamma^2(k)}}\right)^3<\infty$,
	\medskip
	
	\item\label{mthpba:cond4}  $\sum_{t=0}^{\infty}\frac{\gamma(t)\sigma^4(t)}{\sqrt{\sum_{k=t+1}^{\infty}\gamma^2(k)}}<\infty$.
	\medskip
 \end{enumerate}

\revi{Then under Assumptions\s\ref{assum:Ufgradient}-\ref{assum:Ufinfin} and\s\ref{assum:Lip},\ref{assum:infbeh} the sequence $\{\bmu(t)\}$ defined by Algorithm\s\ref{eq:pba2p} converges either to a point from the set of local maxima or to the boundary of one of its connected components almost surely irrespective of the initial parameters.
 Moreover, the agents' joint action updated by the payoff-based procedure defined in Algorithm\s3 converges in probability either to a point from the set of local maxima or to the boundary of one of its connected components.}
\end{thm}
\revi{From the practical point of view, Algorithm 3, similarly to Nesterov's approach in \cite{RePEc:cor:louvco:2011001}, requires agents to coordinate their actions to get the second estimation at the current mean vector $\bmu(t)$. On the other hand, the fact that the variance of the process in Algorithm 3 tends to $0$ under diminishing variance (see \eqref{eq:variance2p} and adjusted definition of $\Wb(t,\zbx(t),\omega)$ in the proof of Theorem\s\ref{th:mthpba}) guarantees a faster convergence of the iterations in comparison with Algorithm 2, where the variance needs to be controlled by the term $\gamma^2(t)$ (see \eqref{eq:finalineq} in the proof of Theorem\s\ref{th:mthpba}).
This conclusion is supported by the simulation results provided in the next section.}

\section{Illustrative Example}\label{sec:impl}
In this section, we design a distributed problem of the following flow control problem \cite{Scutaricdma}. There are $N$ users (agents) in the system, whose goal is to decide on the intensity $a_i\in \R$, $i\in[N]$, of the power flow to be sent over the system. The overall profit in the system is defined by the following function
$$p(a_1,\ldots,a_N)=\log(1 + \sum_{i\in[N]}h_i \exp(a_i))),$$
where $h_i$ corresponds to the reward factor of the flow sent by the user $i$. However, there is some cost
$$c_i(a_i)=3\log(1+\exp(a_i))-a_i$$
the user $i$ needs to pay for choosing the intensity $a_i$, $i\in[N]$. Thus, the objective in the system is to maximize the function
$$\phi(a_1,\ldots,a_N) = p(a_1,\ldots,a_N) - \sum_{i\in[N]}c_i(a_i)$$
over $a_i\in\R$, $i\in[N]$. By using the relation \eqref{eq:pg1} we can conclude that the
problem above, namely
\begin{align}\label{eq:optsim}
 \phi(a_1,\ldots,a_N)\to\max, \quad a_i\in\mathbb R, \quad i\in[N]
\end{align}
can be reformulated in terms of learning potential function maximizers in the potential game $\Gamma=(N, \boldsymbol A=\R^3, \{U_i\}_{i\in[N]}, \phi)$, where
\begin{align*}
 U_i(\ba) = \log\left(1 + \frac{h_i \exp(a_i)}{1+\sum_{j\ne i}h_j \exp(a^j)}\right) - c_i(a_i)
\end{align*}
and $\ba = (a_1,\ldots,a_N)$. In the following we consider the system described above with some positive coefficients $h_i\in (0,1]$, $i\in[N]$.

\subsection{Communication-based approach}
To learn a local maximum of the potential function in $\Gamma$, the users adapt the communi\-cation-based algorithm with the communication topology corresponding to a sequence of random digraphs that is $S$-strongly connected with $S=4$ (see Definition\s\ref{def:CG}).

\begin{figure}[htb!]
\centering
\psfrag{500}[c][b]{\tiny{$500$}}
\psfrag{1500}[c][b]{\tiny{$1500$}}
\psfrag{0}[c][b]{\tiny{$0$}}
\psfrag{1000}[c][b]{\tiny{$1000$}}
\psfrag{2000}[c][b]{\tiny{$2000$}}
\psfrag{2500}[c][b]{\tiny{$2500$}}
\psfrag{3000}[c][b]{\tiny{$3000$}}
\psfrag{4000}[c][b]{\tiny{$4000$}}
\psfrag{3500}[c][b]{\tiny{$3500$}}
\psfrag{-10}[r][c]{\tiny{$-10$}}
\psfrag{-20}[r][c]{\tiny{$-20$}}
\psfrag{-100}[r][c]{\tiny{-100}}
\psfrag{-600}[r][c]{\tiny{$-600$}}
\psfrag{-400}[r][c]{\tiny{$-400$}}
\psfrag{-200}[r][c]{\tiny{$-200$}}
\begin{overpic}[width=0.8\linewidth]{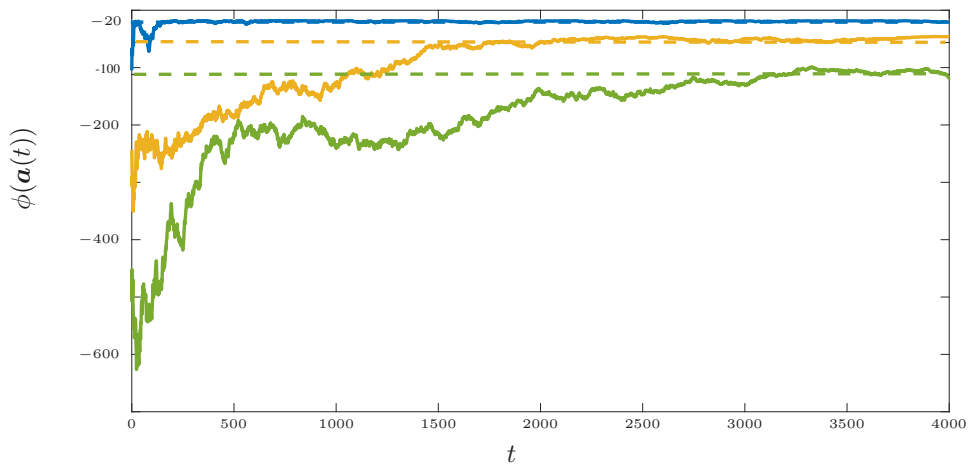}
\put(50.8,-2.4){$t$}
\put(0.8,23){\rotatebox{90}{$\phi(\ba(t))$}}
\end{overpic}
\caption{The value of $\phi$ during the communication-based algorithm for $N=10$ (blue line), $N=25$ (yellow line), and $N=50$ (green line).}
\label{fig:cb200}
\end{figure}

Figures \ref{fig:cb200} demonstrates behavior of the potential function $\phi$ during the run of the communication-based algorithm with some initial users' estimation vectors chosen uniformly at random on the sphere with center at $\boldsymbol 0\in\R^{N}$ and the radius $10$ in the cases $N=10$, $25$, and $50$. As we can see, convergence to a local maximum of the potential function (its value corresponds to the dash lines on the figure) takes place. The convergence rate, however, depends on the number of the agents in the system. The algorithm needs more time to approach a local maximum as $N$ increases.

\subsection{Payoff-based approach}

\begin{figure}[htb!]
\centering
\psfrag{500}[c][b]{\tiny{$500$}}
\psfrag{1500}[c][b]{\tiny{$1500$}}
\psfrag{0}[c][b]{\tiny{$0$}}
\psfrag{1000}[c][b]{\tiny{$1000$}}
\psfrag{2000}[c][b]{\tiny{$2000$}}
\psfrag{2500}[c][b]{\tiny{$2500$}}
\psfrag{3000}[c][b]{\tiny{$3000$}}
\psfrag{4000}[c][b]{\tiny{$4000$}}
\psfrag{3500}[c][b]{\tiny{$3500$}}
\psfrag{-5}[r][c]{\tiny{$-5$}}
\psfrag{-10}[r][c]{\tiny{$-10$}}
\psfrag{-50}[r][c]{\tiny{$-50$}}
\psfrag{-20}[r][c]{\tiny{$-20$}}
\psfrag{-30}[r][c]{\tiny{$-30$}}
\psfrag{-40}[r][c]{\tiny{$-40$}}
\psfrag{-50}[r][c]{\tiny{$-50$}}
\psfrag{-60}[r][c]{\tiny{$-60$}}
\psfrag{-70}[r][c]{\tiny{$-70$}}
\begin{overpic}[width=0.8\linewidth]{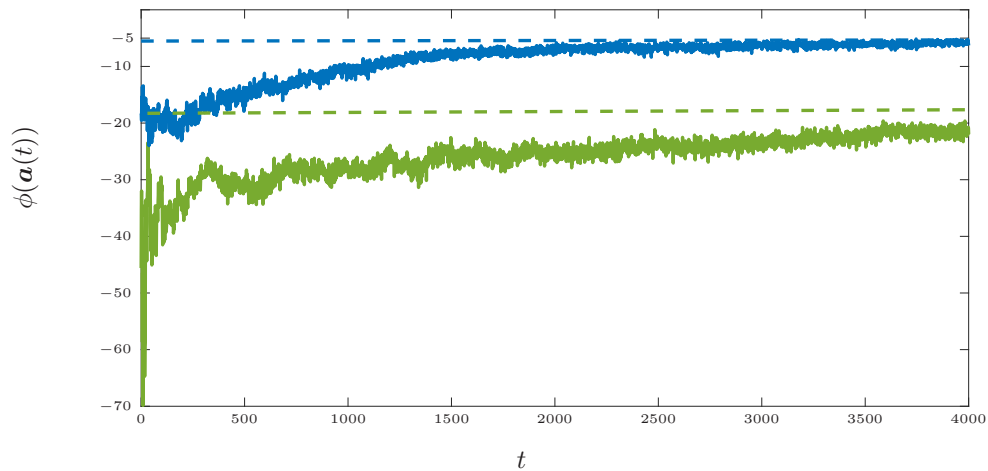}
\put(50.8,-3){$t$}
\put(0,23){\rotatebox{90}{$\phi(\ba(t))$}}
\end{overpic}
\caption{The value of $\phi$ during the payoff-based learning algorithm for $N=3$ (blue line) and $N=5$ (green line).}
\label{fig:pb1000}
\end{figure}

\begin{figure}[htb!]
\centering
\psfrag{500}[c][b]{\tiny{$500$}}
\psfrag{1500}[c][b]{\tiny{$1500$}}
\psfrag{0}[c][b]{\tiny{$0$}}
\psfrag{1000}[c][b]{\tiny{$1000$}}
\psfrag{2000}[c][b]{\tiny{$2000$}}
\psfrag{2500}[c][b]{\tiny{$2500$}}
\psfrag{3000}[c][b]{\tiny{$3000$}}
\psfrag{4000}[c][b]{\tiny{$4000$}}
\psfrag{3500}[c][b]{\tiny{$3500$}}
\psfrag{-150}[r][c]{\tiny{$-150$}}
\psfrag{-50}[r][c]{\tiny{$-50$}}
\psfrag{-100}[r][c]{\tiny{$-100$}}
\psfrag{-300}[r][c]{\tiny{$-300$}}
\psfrag{-250}[r][c]{\tiny{$-250$}}
\psfrag{-200}[r][c]{\tiny{$-200$}}
\begin{overpic}[width=0.8\linewidth]{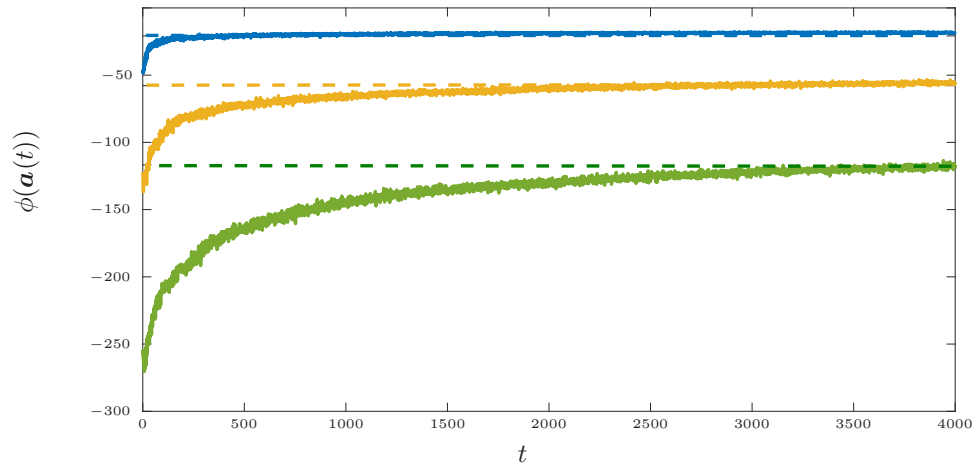}
\put(50.8,-2){$t$}
\put(0.0,23){\rotatebox{90}{$\phi(\ba(t))$}}
\end{overpic}
\caption{The value of $\phi$ during the payoff-based learning algorithm for $N=10$ (blue line), $N=25$ (yellow line), and $N=50$ (green line).}
\label{fig:pb200}
\end{figure}

If communication does not take place and the gradient information is not available, the agents use the payoff-based algorithm.
Figures\s\ref{fig:pb1000} and \ref{fig:pb200} demonstrate the performance of the algorithm based on one- and two-point estimations respectively, given some initial mean vector $\bmu(0)$ chosen uniformly at random on the sphere with the center at $\boldsymbol 0\in\R^{N}$ and the radius $10$.
The simulations support the discussion in Subsection\s\ref{subsec:pb2p}. Namely, the payoff-based algorithm based on one estimation of the utilities requires a more sophisticated tuning of the parameter $\sigma(t)$ and $\gamma(t)$ to approach a local maximum within some finite time.
As we can see on Figure\s\ref{fig:pb1000}, this algorithm needs over $4000$ iterations to get close to a local maximum (the value of the potential function at the local maximum corresponds to the dash line) already in the case $N=5$. To overcome this practical limitation we use the procedure using two estimations of the utilities' values per iteration (see Subsection\s\ref{subsec:pb2p}). The corresponding simulation results are presented on Figure\s\ref{fig:pb200} for the cases $N=10$, $25$, and $50$. As we can see, the procedure converges to a local maximum within some finite time (as before, the dash lines demonstrate the values of the potential function at the local maxima).

\section{Conclusion}\label{sec:concl}
In this paper we introduced the stochastic learning algorithms that can be applied to optimization problems modeled by potential games with continuous
actions. We focused on two different types of information. The first learning procedure uses the non-restrictive push-sum protocol to enable communication between agents. The second algorithm is payoff-based one. The advantage of this procedure is that players need neither communication, nor knowledge on the analytical properties of their utilities to follow
the rules of the algorithm. Without any assumption on concavity of the potential function, we demonstrated that under appropriate settings of the algorithms, agents' joint actions converge to a critical point different from local minima. \revi{Future work can be devoted to investigation of escaping from saddle points in distributed multi-agent optimization problems with different information restrictions in systems as well as rigorous analysis of the convergence rates.}
\bibliographystyle{plain}
\bibliography{format2}

\end{document}